%% file: submittedetds.tex
\documentclass{etds}    

\usepackage{amsfonts}
\usepackage{amssymb}
\usepackage{amsmath}
\usepackage{enumerate}

\usepackage[latin2]{inputenc}
\usepackage{graphicx}

\def\eps{\varepsilon}
\def\phi{\varphi}

\def\N{\mathbb{N}}
\def\Q{\mathbb{Q}}

\def\phik{\phi_1, \phi_2, \ldots, \phi_r}
\def\psik{\psi_1, \psi_2, \ldots, \psi_s}
\def\={\overset{\textrm{def}}{=}}
\def\lep{\;\dot{\le}\;}
\def\gep{\;\dot{\ge}\;}
\def\bigcupp{\sideset{}{^*}\bigcup}
\def\bigcuppp{\sideset{}{^{**}}\bigcup}
\def\cupp{\cup^{*}}
\def\cuppp{\cup^{**}}
\def\affk{A_K}  

\def\egyenlet{K=\phi_1(K) \cup \ldots \cup \phi_r(K)}
\def\egy{K=\phi_1(K) \cupp \ldots \cupp \phi_r(K)}
\def\egys{K=\psi_1(K) \cupp \ldots \cupp \psi_s(K)}
\def\phiIi{\phi_{I_i}}
\def\phiIK{\phi_I(K)}

\def\phiJi{\phi_{J_i}}
\def\phiIiK{\phi_{I_i}(K)}

\def\ubolv{\overline{U}\setminus V}
\def\minid{\setminus\{\textrm{identity}\}}
\def\bol{\setminus}
\def\id{\{\textrm{identity}\}}
\def\ii{\textrm{identity}}
\def\dist{\textrm{dist}}
\def\alfa{\alpha}
\def\phii{\phi^{-1}}
\def\ures{\emptyset}
\def\noteq{\neq}
\def\muH{\mu_H}
\def\G{\mathcal{G}}
\def\F{\mathcal{F}}

\def\Ak{\mathcal{A}_K}
\def\S{\mathcal{S}}
\def\Sk{\mathcal{S}_K}

\def\Ik{\mathcal{I}_K}
\def\l{\lambda}
\def\lm{\l_m}

\def\kb{\approx}
\def\int{\textrm{int}}

\providecommand{\mult}[1]{\begin{multline} #1 \nonumber \end{multline}}

\providecommand{\M}[1]{\mu\Big(#1\Big)}
\providecommand{\m}[1]{\mu\big(#1\big)}
\providecommand{\mh}[1]{\mu_H\big(#1\big)}
\providecommand{\mmcs}[1]{\mu^{*}(#1)}
\providecommand{\mm}[1]{\mu(#1)}

\providecommand{\mcs}[1]{\mu^{*} \big(#1\big)}

\newtheorem{theorem}{Theorem}[section]
\newtheorem{llll}[theorem]{Lemma}
\newtheorem{lemma}[theorem]{Lemma}
\newtheorem{facts}[theorem]{Facts}
\newtheorem{propo}[theorem]{Proposition}
\newtheorem{ccc}[theorem]{Claim}
\newtheorem{conjecture}[theorem]{Conjecture}
\newtheorem{eee}[theorem]{Example}
\newtheorem{rrr}[theorem]{Remark}
\newtheorem{sss}[theorem]{Statement}
\newtheorem{question}[theorem]{Question}
\newtheorem{cccc}[theorem]{Corollary}
\newtheorem{ddd}[theorem]{Definition}

\newtheorem{notation}[theorem]{Notation}
\newtheorem{tetel}[theorem]{Theorem}
\newtheorem{all}[theorem]{Proposition}
\newtheorem{kov}[theorem]{Corollary}

\newcommand{\sumi}{\sum_{i=1}^\infty}
\newcommand{\cupi}{\cup_{i=1}^{\infty}}
\newcommand{\cupsi}{\cup_{i=1}^{*\infty}}
\newcommand{\sumj}{\sum_{j=1}^\infty}

\newcommand{\cupsj}{\cup_{j=1}^{*\infty}}
\newcommand{\norm}{||}
\providecommand{\norma}[1]{\lVert #1 \rVert} 
\newcommand{\sm}{\setminus}

\newcommand{\bt}{\begin{theorem}}
\newcommand{\bn}{\begin{notation}\upshape}
\newcommand{\bl}{\begin{llll}}
\newcommand{\bc}{\begin{ccc}}
\newcommand{\bex}{\begin{eee}}
\newcommand{\br}{\begin{rrr}\upshape}

\newcommand{\bs}{\begin{sss}}
\newcommand{\bd}{\begin{ddd}\upshape}
\newcommand{\bq}{\begin{qqq}}
\newcommand{\bcoro}{\begin{cccc}}

\newcommand{\bp}{\proc{Proof.}}

\def\Biz{\bp}
\def\qed{\ep\medbreak}

\newcommand{\et}{\end{theorem}}
\newcommand{\en}{\end{notation}}
\newcommand{\el}{\end{llll}}
\newcommand{\ec}{\end{ccc}}
\newcommand{\eex}{\end{eee}}
\newcommand{\er}{\end{rrr}}
\newcommand{\es}{\end{sss}}
\newcommand{\ed}{\end{ddd}}
\newcommand{\eq}{\end{qqq}}
\newcommand{\ecoro}{\end{cccc}}

\newcommand{\lab}[1]{\label{#1}}

\newcommand{\RR}{\mathbb{R}}
\newcommand{\R}{\mathbb{R}}
\newcommand{\Rd}{\mathbb{R}^d}

\newcommand{\Z}{\mathbb{Z}}

\newcommand{\su}{\subset}

\newcommand{\be}{\beta}

\newcommand{\al}{\alpha}

\newcommand{\de}{\delta}

\newcommand{\si}{\sigma}

\newcommand{\la}{\lambda}
\newcommand{\mut}{\tilde{\mu}}

\newcommand{\iM}{\mathcal{M}}
\newcommand{\MA}{\mathcal{M}_A}

\newcommand{\Rn}{\RR^n}

\newcommand{\diam}{\mathrm{diam}} 

\begin{document}
\ETDS{0}{0}{0}{0}

\runningheads{M.~Elekes, T.~Keleti, A.~M\'ath\'e}{Self-similar and self-affine sets}

\title{Self-similar and self-affine sets; 
measure of the intersection of two copies}

\author{{M\'arton Elekes\affil{1}}\footnote{Supported by Hungarian Scientific Foundation grant no.~37758.},
{Tam\'as Keleti\affil{2}}\footnote{Supported by Hungarian Scientific Foundation grant no.~F 43620.}
\ and
{Andr\'as M\'ath\'e\affil{2}}\footnote{Supported by Hungarian Scientific Foundation grant no.~T 49786.}
}

\address{\affilnum{1}\ Alfr\'ed R\'enyi Institute of Mathematics, Hungarian Academy of
Sciences,
P.O.~Box 127, H-1364, Budapest, Hungary\\
\affilnum{2}\ Department of Analysis, E\"otv\"os Lor\'and University,
P\'az\-m\'any P\'e\-ter s\'et\'any 1/c, H-1117 Budapest, Hungary\\
\email{emarci@renyi.hu, elek@cs.elte.hu, amathe@cs.elte.hu}}

\recd{2008}
\begin{abstract}
Let $K\su\Rd$ be a self-similar or self-affine set, let $\mu$ be a self-similar
or self-affine measure on it, and let $\G$ be the group of affine maps,
similitudes, isometries or translations of $\Rd$. Under various assumptions
(such as separation conditions or we assume that the transformations
are small perturbations or that $K$ is a so called Sierpi\'nski sponge) we
prove theorems of the following types, which are closely related to each other;
\begin{itemize}
\item \emph{(Non-stability)}

\noindent
There exists a constant $c<1$ such that for every $g\in\G$ we have either
$\m{K\cap g(K)}<c\cdot\mu(K)$
or $K\su g(K)$.

\item \emph{(Measure and topology)}

\noindent 
For every $g\in\G$ we have $\m{K\cap g(K)} > 0 \iff \int_K (K\cap g(K)) \neq
\emptyset$ (where $\int_K$ is interior relative to $K$). 

\item \emph{(Extension)}

\noindent 
The measure $\mu$ has a $\G$-invariant extension to $\R^n$.
\end{itemize}

\noindent Moreover, in many situations we characterize those $g$'s for which
$\m{K\cap g(K)} > 0$ holds, 
and we also get results about those $g$'s for which $g(K)\su K$ or 
$g(K)\supset K$ holds.

\end{abstract}






\section{Introduction}
The study of the size of the intersection of Cantor sets
has been a central research area in geometric measure theory
and dynamical systems lately, see e.g. the works
of Igudesman \cite{I}, Li and Xiao \cite{LX}, Moreira \cite{Mo},
Moreira and Yoccoz \cite{MY}, Nekka and Li \cite{NL}, 
Peres and Solomyak \cite{PS}.
For instance J-C. Yoccoz and C. G. T. de Moreira \cite{MY} proved
that if the sum of the Hausdorff dimensions of two regular Cantor sets
exceeds one then, in the typical case, there are translations of them 
stably having intersection with positive Hausdorff dimension.

The main purpose of this paper is to study the measure of the intersection of
two Cantor sets which are (affine, similar, isometric or translated)
copies of a self-similar or self-affine set in $\Rd$.
By measure here we mean a self-similar or self-affine measure on one of
the two sets.

We get instability results 
stating that the measure of the intersection is separated from
the measure of one copy.
This strong non-continuity property is in sharp contrast with the 
well known fact that for any Lebesgue measurable set $H\su\Rd$ 
with finite measure the 
Lebesgue measure of $H\cap (H+t)$ is continuous in $t$.

We get results stating that the intersection is of
positive measure if and only if it contains a relative open set.
This result resembles some recent deep results 
(e.g. in \cite{LW}, \cite{MY})
stating that for certain classes of sets having positive Lebesgue measure
and nonempty interior is equivalent.
In the special case when the self-similar set is the classical Cantor set
our above mentioned results were obtained by F. Nekka and Jun Li \cite{NL}.
For other related results see also the work of
Falconer \cite{Fa}, Feng and Wang \cite{FW}, 
Furstenberg \cite{Fu}, Hutchinson \cite{H}, 
J\"arvenp\"a\"a \cite{J} and Mattila \cite{M1}, \cite{M2}, \cite{M3}.

As an application we also get isometry (or at least translation)
invariant measures of $\R^d$ such that the measure of the given self-similar
or self-affine set is $1$.

Feng and Wang \cite{FW} has proved recently ``The Logarithmic Commensurability
Theorem'' about the similarity ratios of a homogeneous self-similar set in
$\R$ with the open set condition and a similarity map that maps the
self-similar set into itself 
(see more precisely after Theorem~\ref{lambdatul}), and they also posed the
problem of generalizing their result to higher dimensions. 
For self-similar sets with the strong separation condition we prove a higher
dimensional generalization without assuming homogeneity.

\subsection{Self-affine sets.}

Let $K\su\Rd$ be a \emph{self-affine set} with the 
\emph{strong separation condition};
that is, $K=\phi_1(K)\cupp\ldots\cupp \phi_r(K)$ is a compact set, where 
$r\ge 2$ and $\phi_1,\ldots,\phi_r$
are injective and contractive $\Rd\to\Rd$ affine maps and $\cupp$ denotes
disjoint union.

For any $p_1,\ldots,p_r\in(0,1)$ such that $p_1+\ldots+p_r=1$ 
let $\mu$ be the corresponding 
\emph{self-affine measure}; that is,
the image of the 
infinite product of  the discrete probability measure $p(\{i\})=p_i$ 
on $\{1,\ldots,r\}$
under the representation map
$
\pi:\{1,\ldots,r\}^{\N}\to K,\quad 
\{\pi(i_1,i_2,\ldots)\}=\cap_{n=1}^\infty 
(\phi_{i_1}\circ\ldots\circ\phi_{i_n})(K)
$.

In Section~\ref{affine} we show (Theorem~\ref{kovesarki}) that 
small affine perturbations of $K$ cannot intersect a very large part of $K$;
that is, there exists a $c<1$ and 
a neighborhood $U$ of the identity map in the space
of affine maps such that for any $g\in U\sm\id$ 
we have $\m{K\cap g(K)}<c$.
We also prove (Theorem~\ref{selfaffine}) that no isometric but nonidentical 
copy of $K$ can intersect a very large part of $K$; that is, 
there exists a constant $c<1$ such that for any isometry $g$ either
$\m{K\cap g(K)}<c$ or $g(K)=K$.


\subsection{Self-similar sets.}

Now let $K\su\Rd$ be a \emph{self-similar set} with the 
\emph{strong separation condition} and $\mu$ a self-similar measure
on it; that is, $K$ and $\mu$ are defined as above 
with the extra assumption that
$\phi_1,\ldots, \phi_r$ are similitudes.

In Section~\ref{similar}  we prove (Theorem~\ref{andras})
that for any given self-similar set $K\su\R^d$ with the
strong separation condition and self-similar measure $\mu$ on $K$ 
there exists a $c<1$ such that for any similitude $g$
either $\m{K\cap g(K)}<c\cdot\mu(K)=c$ or $K\su g(K)$. In other words, the
intersection of a self-similar set with the strong separation condition
and its similar copy cannot have a really big non-trivial intersection.

Let $K$, $\mu$ and $g$ be as above. 
An obvious way of getting 
$\m{K\cap g(K)}>0$ is when $g(K)$ contains a nonempty
(relative) open set in $K$. 
The main result (Theorem~\ref{pos}) of Section~\ref{similar}, 
which will follow from the above mentioned Theorem~\ref{andras}, shows that 
this is the only way.
That is, for any self-similar set $K\su\R^d$ with the
strong separation condition and self-similar measure $\mu$ on $K$
a similar copy of $K$ has positive $\mu$ measure in $K$ if and only if
it has nonempty relative interior in $K$.

An immediate consequence (Corollary~\ref{mindegy}) of the 
above result is that for any fixed self-similar set with the
strong separation condition and for any two self-similar measures
$\mu_1$ and $\mu_2$ we have 
$\mu_1\big(g(K)\cap K\big)>0 \Longleftrightarrow\mu_2\big(g(K)\cap K\big)>0$
for any similitude $g$.
As an other corollary (Corollary~\ref{countable})
we
get that for any given self-similar set $K\su\Rd$ with the strong
separation condition and self-similar measure $\mu$ on $K$ there exist
only countably many (in fact exactly countably infinitely many)
similitudes $g:\affk\to\R^d$ 
(where $\affk$ is the affine span of $K$)
such that $g(K)\cap K$ has positive $\mu$-measure.

Let $K\su\Rd$ be a self-similar set with the strong
separation condition and let $s$ be its  
Hausdorff dimension, which in this case equals its
similarity and box-counting dimension.
Then the $s$-dimensional Hausdorff measure is a constant multiple
of a self-similar measure (one has to choose $p_i=a_i^s$, where
$a_i$ is the similarity ratio of $\phi_i$).
Therefore all the above results hold 
when $\mu$ is $s$-dimensional Hausdorff measure.

In Section~\ref{similar} we also need and get results
(Proposition~\ref{haus}, Lemma~\ref{vegessok}, 
Theorem~\ref{lambdatul} and Corollary~\ref{log})
stating that only very special similarity maps can map a self-similar set with
the strong separation condition into itself. 
Theorem~\ref{lambdatul} and Corollary~\ref{log} are the already mentioned
generalizations of The Logarithmic Commensurability Theorem of 
Feng and Wang \cite{FW}.
 
In Section~\ref{measures} we apply the main result  
(Theorem~\ref{pos}) and some of the above mentioned results
(Lemma~\ref{vegessok} and Theorem~\ref{lambdatul}) 
of Section~\ref{similar} to characterize those 
self-similar measures on a self-similar set with the strong 
separation condition that can be extended to $\Rd$ as an isometry
invariant Borel measure. 
It turns out that, unless there is a clear obstacle,
any self-similar measure can be extended to $\Rd$ as an isometry
invariant measure. 
Thus, for a given self-similar set with the strong 
separation condition,
there are  usually many distinct isometry invariant Borel measures
for which the set is of measure $1$.

Let us simply call a measure defined on $K$ isometry invariant if it can be extended to an isometry invariant measure on $\R^d$.
Many different collections of similitudes can define the same self-similar set. 
We call $\{\phik\}$ a \emph{presentation} of $K$ if $\egy$ holds; in other
words, $K$ is the attractor of the iterated function system $\{\phik\}$ with the extra condition of disjointness.

The notion of a self-similar measure
on $K$ depends on the particular presentation. However, we show that the notion of
isometry invariant self-similar measure on $K$ is indifferent of the presentations
(Theorem~\ref{vanertelme}).
By this theorem we
can define a natural number for each self-similar set (satisfying 
the strong separation property), an invariant,
which does not depend on the presentation (Theorem~\ref{algebraic}). This invariant is equal to the dimension of the space of isometry invariant self-similar measures, and is related to the algebraic dependence of the similitudes of some (any) presentation of $K$.

In Section~\ref{pelda} we show that the connection between different presentations of a self-similar set can be very complicated. 
This sheds some light on why results and their proofs 
in Section~\ref{measures} are complicated.
The structure of different presentations of a self-similar set in $\R$ has been
also studied recently and independently by Feng and Wang in \cite{FW}, where
a similar example is presented.

%
\subsection{Self-affine sponges.}
Take the $[0,1]^n$ unit cube in $\R^n$ ($n\in\N$) and subdivide it into 
$m_1\times\ldots\times m_n$ boxes of same size 
($m_1,\ldots,m_n\ge 2$) and cut out some of them. 
Then do the same with the remaining boxes using the same
pattern as in the first step and so on. What remains after infinitely
many steps is a self-affine set, which is called
\emph{self-affine Sierpi\'nski sponge}. (A more precise definition will be given
in Definition~\ref{def:K}.)

For $n=2$ these sets were studied in several papers (in which they
were called self-affine carpets or self-affine carpets of Bedford
and McMullen). Bedford \cite{Be} and McMullen \cite{Mu} determined
the Hausdorff and Minkowski dimensions of these self-affine carpets.
(The Hausdorff and Minkowski
dimension of self-affine Sierpi\'nski sponges was determined
by Kenyon and Peres \cite{KP}).
Gatzouras and Lalley \cite{GL} proved that except in some 
relatively simple
cases such a set has zero or infinity Hausdorff measure in its
dimension (and so in any dimension). Peres extended their results
by proving that (except in the same rare simple cases) for any gauge
function neither the Hausdorff \cite{PeH} nor the packing \cite{PeP}
measure of a self-affine carpet can be positive and finite
(in fact, the packing measure cannot be $\si$-finite either),
and remarked that these results extend to self-affine 
Sierpi\'nski sponges of higher dimensions.

Recently the first and the second listed authors of the present paper
showed \cite{EK} that 
some nice sets -- among others the set of Liouville numbers --
have zero or non-$\sigma$-finite Hausdorff and packing measure
for any gauge function by proving that these sets have zero 
or non-$\sigma$-finite measure for any translation invariant Borel 
measure. (Much earlier Davies \cite{Da} constructed a compact 
subset of $\R$ with this property.) So it was natural to ask
whether the self-affine carpets of Bedford and McMullen have this
stronger property.

In Section~\ref{intersection} we prove (Corollary~\ref{interior}) that 
for any self-affine 
Sierpi\'nski sponge 
$K\su\R^n$ ($n\in\N$) with the natural Borel
probability measure $\mu$ (see in Definition~\ref{standard})
on $K$ and $t\in\R^n$,
the set $K\cap(K+t)$ has positive $\mu$ measure
if and only if it has non-empty interior relative to $K$.

For this we prove (Theorem~\ref{structure})
that for any self-affine Sierpi\'nski sponge
$K\su\R^n$ ($n\in\N$) and translation vector $t\in\R^n$ we
have $\m{K\cap(K+t)}=0$ 
unless $K$ or $t$ are of very special form.

We also characterize (Theorem~\ref{instabilsponge}) those 
Sierpi\'nski sponges for which we do not have instability result
for translations and the natural probability measure $\mu$. In fact, we get that 
$\m{K\cap(K+t)}$ can be close to $1$ only for the same special
sponges that appear in the above mentioned result.

In Section~\ref{spongemeasure} we show (Theorem~\ref{main}) that
for any self-affine Sierpi\'nski sponge $K\su\Rn$ 
the natural probability measure $\mu$ on $K$ can be extended as a
translation invariant Borel measure $\nu$ on $\R^n$.
We also extend this result (Theorem~\ref{direct}, Corollary~\ref{cosc})
to slightly larger classes of self-affine sets.

\section{Notation, basic facts and some lemmas}

In this section we collect several notions and well known or fairly easy 
statements that we will need in the sequel.
Some of these might be interesting in their own right. 
Of course, only a few of them are needed for each specific section.
Though some of these statements may be well known, for the sake
of completeness we included the proofs.

\begin{notation}\upshape
We shall denote by $\cupp$ the disjoint union
and by $\dist$ the (Euclidean) distance.
\end{notation}

\subsection{Affine maps, similitudes, isometries.}

\bd
A mapping $g:\Rd\to\Rd$ is called a \emph{similitude} 
if there is a constant $r>0$,
called \emph{similarity ratio}, such that 
$\dist(g(a),g(b))=r\cdot \dist(a,b)$ for any $a,b\in\Rd$.

The \emph{affine map}s of $\Rd$ are of the form 
$x \mapsto A x + b$, where 
$A$ is an $n\times n$ matrix and
$b\in\R^d$ is a translation vector. 
Thus the set of all affine maps
of $\Rd$ can be considered as $\R^{d^2+d}$ and
so it can be considered as a metric space. 
\ed

It is easy to check that a sequence 
$(g_n)$ in this metric space
converges to an affine map $g$ if and only if $g_n$ converges to
$g$ uniformly on any compact subset of $\Rd$.

\bd
For a given set $K\su\Rd$ with affine span $\affk$ let $\Ak$, $\Sk$ and $\Ik$ 
denote the metric space (with the above metric) 
of the injective affine maps, similitudes and
isometries of $\affk$, respectively.
\ed

Note also that all these three metric spaces with the composition
can be also considered as topological groups.

\subsection{Self-similar and self-affine sets and measures.}

\bd
A $K\su\Rd$ compact set is a \emph{self-similar/self-affine set} 
if $K=\phi_1(K)\cup\ldots\cup \phi_r(K)$, where 
$r\ge 2$ and $\phi_1,\ldots,\phi_r$
are similitudes/injective and contractive affine maps.

By the \emph{$n$-th generation elementary pieces of $K$} we mean the
sets of the form $(\phi_{i_1}\circ\ldots\circ\phi_{i_n})(K)$, 
where $n=0,1,2,\ldots$.

We shall use multi-indices. 
By a \emph{multi-index} we mean a finite sequence
of indices; for $I=(i_1, i_2, \ldots, i_n)$ let
$\phi_I = \phi_{i_1}\circ\ldots\circ\phi_{i_n}$ 
and $p_I = p_{i_1}p_{i_2}\ldots p_{i_n}$.
We shall consider $I=\emptyset$ as a multi-index as well:
$\phi_\emptyset$ is the identity map and $p_\emptyset = 1$.
\ed

Note that the elementary pieces of $K$ are
the sets of the form $\phi_I(K)$. These sets are also 
self-similar/self-affine;
and if $h$ is an injective affine map then $h(K)$ is also 
self-similar/self-affine and
its elementary pieces are the sets of the form $h(\phi_I(K))$.


\bd\label{defmeasure}
Let $\egyenlet$ be a self-similar/self-affine set, and let
$p_1+\ldots+p_r=1$, $p_i>0$ for all $i$.
Consider the symbol space $\Omega=\{1,\ldots,r\}^\N$ 
equipped with the product topology and 
let $\nu$ be the Borel measure on $\Omega$ 
which is the countable infinite product of the 
discrete probability measure $p(\{i\})=p_i$ on $\{1,\ldots,r\}$.
Let $$
\pi:\Omega\to K,\quad 
\{\pi(i_1,i_2,\ldots)\}=\cap_{n=1}^\infty 
(\phi_{i_1}\circ\ldots\circ\phi_{i_n})(K)
$$
be the continuous addressing map of $K$.
Let $\mu$ be the image measure of $\nu$ under the projection $\pi$; that is,
\begin{equation}\label{altdef}
\mu(H)=\nu\big(\pi^{-1}(H)\big) \quad \textrm{for every Borel set } H\subset K.
\end{equation}

Such a $\mu$ is called a 
\emph{self-similar/self-affine measure}
on $K$.
\ed

One can also define (see e.g. in \cite{Fa97}) self-similar or self-affine 
measures as the
unique probability measure $\mu$ on $K$ such that
$$\mu(H)=\sum_{i=1}^r p_i \m{\phi_i^{-1}(H)}$$
holds for every Borel set $H\su K$.
It was already proved by Hutchinson \cite{H} that the two definition agrees.

\begin{lemma}\label{jancski}
Let $\egyenlet$ be a self-affine set, $p_1+\ldots+p_r=1$, $p_i>0$ for all $i$, 
and let $\mu$ be the self-affine measure on $K$ 
corresponding to the weights $p_i$.

Then for every affine subspace $A$ either 
$\mu(A\cap K)=0$ or $A\supset K$.
\end{lemma}

\bp
Let $\{x_1, x_2, \ldots, x_k\}$ 
be a maximal collection of affine independent
points in $K$. Choose $U_1,\ldots,U_k$ convex open 
sets such that $x_j\in U_j$ ($j=1,\ldots,k$) and whenever
we choose one point from each $U_j$ they are affine independent.
Since $K\cap U_i$ is a nonempty relative open subset of $K$,
we may choose an elementary piece $\phi_{I_j}(K)$ in $U_j$ for each $j$.
Let $\eps= \min_{1\le j\le k} p_{I_j}>0$.

We shall use the notation we introduced in Definition~\ref{defmeasure}.
For $1\le i\le r$ and $\omega=(i_0, i_1, \ldots)\in \Omega$, let $\sigma_i(\omega)=(i, i_0, i_1, \ldots)$.
Thus $\nu\big(\sigma_i(H)\big)=p_i \nu(H)$ for all Borel subset $H$ of $\Omega$.

Suppose that $A$ is an affine subspace such that $\mu(A\cap K)>0$. 
Thus $\nu\big(\pi^{-1}(A)\big)>0$. 
It is easy to prove 
(see a possible argument later in the proof of Lemma~\ref{suruseg})
that this implies that 
there exists an elementary piece $\sigma_J(\Omega)$ such that
$$
\nu\big(\pi^{-1}(A) \cap \sigma_{J}(\Omega)\big) > (1-\eps) \nu\big(\sigma_{J}(\Omega)\big)= (1-\eps) p_J.
$$
Since $\nu\big((\sigma_J\circ\sigma_{I_j})(\Omega)\big)=p_J p_{I_j} \ge p_J \eps$ ($j=1,\ldots, k$),
the set $\pi^{-1}(A)$ must intersect the sets $(\sigma_J\circ\sigma_{I_j})(\Omega)$.
Therefore the set $A$ must intersect the sets $\pi((\sigma_J\circ\sigma_{I_j})(\Omega))=(\phi_J\circ\phi_{I_j})(K)$
($j=1,\ldots, k$).

By picking one point from each $A\cap(\phi_J \circ \phi_{I_j})(K)$,
we get a maximal collection of affine independent points 
in $K$ since $\phi_J$ is an invertible affine mapping. 
As this collection is contained in the affine subspace $A$,
we get that $K$ is also contained in $A$.
\qed

\br
In this paper one of our main goals is to study  
$\m{K\cap g(K)}$, where $g$ is an affine map of $\Rd$.
By the above lemma if the affine map $g$ does not map
the affine span $\affk$ of $K$ onto itself then 
$\m{g(K)\cap K}=0$ since $K\not\subset g(\affk)$.
The other property of affine maps we are interested in is $K\su g(K)$, which
also implies that $g$ maps $\affk$ onto itself.
Thus it is enough to consider those affine maps $g$ of $\Rd$
that map the affine span $\affk$ of $K$ onto itself.
Since then both $K$ and $g(K)$ are in $\affk$, only the restriction 
of $g$ to $\affk$ matters. 
This is why in the next section we shall study $\Ak$, $\Sk$ and $\Ik$ 
(the injective affine maps, similitudes and isometries of $\affk$) 
instead of all affine maps, similitudes 
and isometries of $\Rd$.

Therefore if we state something (about $\m{g(K)\cap K}$ or
about the property $K\su g(K)$) for every affine map, similitude
or isometry $g$,
it will be enough to prove them for $g\in\Ak$, $g\in\Sk$ or $g\in\Ik$, 
respectively.

Note also that self-similar sets and measures are self-affine as well,
so results about self-affine sets and measures also apply for
self-similar sets and measures.
\er

\subsection{Separation properties.}

\bd
A self-similar/self-affine set $\egyenlet$ (or more precisely, the
collection $\phi_1,\ldots,\phi_r$ of the representing maps) satisfies the 
\begin{itemize}
\item \emph{strong separation condition} (SSC) if the union 
$\phi_1(K)\cupp \ldots \cupp \phi_r(K)$ is disjoint;
\item \emph{open set condition} (OSC) if there exists a nonempty 
bounded open set $U\su\Rd$ such that
$\phi_1(U)\cupp \ldots \cupp \phi_r(U)\su U$;
\item \emph{strong open set condition} (SOSC) if there exists a nonempty 
bounded open set $U\su\Rd$ such that $U\cap K\neq \emptyset$ and
$\phi_1(U)\cupp \ldots \cupp \phi_r(U)\su U$;
\item \emph{convex open set condition} (COSC) if there exists a nonempty

bounded open convex set $U\su\Rd$ such that
$\phi_1(U)\cupp \ldots \cupp \phi_r(U)\su U$;
\item \emph{measure separation condition} (MSC) if for any 
self-similar/self-affine measure $\mu$ on $K$ we have 
$\m{\phi_i(K)\cap \phi_j(K)}=0$ for any $1\le i<j\le r$.
\end{itemize}
\ed

We note that the first three definitions are standard but we have not
seen any name for the last two in the literature.

It is easy to check the well known fact that we must have $K\su \overline{U}$
(where $\overline E$ denotes the closure of a set $E$) for the open
set $U$ in the definition of OSC (and SOSC, COSC).

It is easy to see ($U$ can be chosen as a small $\eps$-neighborhood
of $K$ for the first implication) that for any self-affine set
$$
SSC \Longrightarrow SOSC \Longrightarrow OSC.
$$
Using the methods of C. Bandt and S. Graf \cite{BG}, 
A. Schief proved in \cite{S} that, 
in fact, $SOSC \Longleftrightarrow OSC$ holds
for self-similar sets.

In \cite{S} for self-similar sets $SOSC \Longrightarrow MSC$ is also proved.
Since the proof works for self-affine sets as well we get that
for any self-affine set
$$
SOSC \Longrightarrow MSC.
$$

It seems to be also true that $COSC \Longrightarrow SOSC$ and so
$COSC \Longrightarrow MSC$ but we do not prove this, since we do not
need the first implication and the following lemma is stronger than the
second implication.

\bl\label{cosc-msc}
Let $\egyenlet$ be a self-affine set in $\R^d$ 
with the convex open set condition
and let $\mu$ be a self-affine measure on it. 
Then for any affine map $\Psi:\Rd\to\Rd$ we have
$$
\mu\Big(\Psi\big( \phi_i(K)\cap \phi_j(K) \big) \Big)=0 \qquad (\forall \,1\le i < j \le r).
$$
\el

\bp 
Let $1\le i < j \le r$ and $U$ be the convex open set given in the 
definition of COSC. 
Let $\affk$ be the affine span of $K$.
Since $\phi_i(U\cap \affk)$ and $\phi_j(U\cap \affk)$ are disjoint convex
open sets in $\affk$, 
$\overline{\phi_i(U\cap \affk)}\cap\overline{\phi_j(U\cap \affk)}$
must be contained in a proper affine subspace $A$ of $\affk$. 
Since $K\su\overline{U}\cap \affk$, this implies that 
$\phi_i(K)\cap \phi_j(K)\su A$, and so
\begin{equation}\label{PsiA}
\Psi\big( \phi_i(K)\cap \phi_j(K) \big)\su \Psi(A).
\end{equation}
Since $\Psi(A)$ is an affine subspace, 
which is smaller dimensional than the affine span $\affk$ of $K$,
we cannot have $K\su \Psi(A)$, so by Lemma~\ref{jancski} we must
have $\m{K\cap \Psi(A)}=0$. By (\ref{PsiA}) this implies that
$\m{\Psi( \phi_i(K)\cap \phi_j(K) ) }=0$.
\qed

We also note that one can 
find a self-similar set in $\R$
that satisfies even the SSC but does not satisfy the COSC \cite[Example 5.1]{FW}, 
so SSC and COSC are independent even for self-similar sets of $\R$.

\bn
Given a fixed measure $\mu$, we shall say that
two sets are \emph{almost disjoint} if their intersection
has $\mu$-measure $0$. 
The almost disjoint union will be denoted by $\cuppp$.
\en

It is very easy to prove one by one each of the following facts.

\begin{facts}\label{facts}
Let $\egyenlet$ be a self-affine/self-similar set with the measure separation 
condition and let $\mu$ be a self-affine/self-similar measure on it,
which corresponds to the weights $p_1,\dots,p_r$. Then the 
following statements hold.
\begin{enumerate}
\item Any two elementary pieces of $K$ are either almost disjoint or
one contains the other.
\item Any union of elementary pieces can be replaced by an almost disjoint countable union.
\item For any multi-index $I$ we have $\mu\circ\phi_I=p_I\cdot \mu$;
that is, $\mu\circ\phi_I(B)=p_I\cdot \mu(B)$ for any Borel set $B\su K$. 
\item We have $\m{\phi_I(K)}=p_I$ for any multi-index $I$.
\item \label{origdef} For any Borel set $B\su K$ we have
$$
\mu(B)=\inf\Big\{\sum_{i=1}^\infty p_{I_i}: 
B \subset \bigcuppp_{i=1}^{\infty} \phi_{I_i}(K)\Big\}.  
$$
\end{enumerate}
\end{facts}





Since SOSC and COSC are both stronger than MSC and one of them 
will be always assumed in this paper, the statements of this
lemma will often be tacitly used.
Sometimes, for example, we shall even handle 
the above almost disjoint sets as disjoint sets and often consider
Fact~\ref{origdef} as the definition of self-affine/self-similar measures.

\begin{lemma}\label{suruseg}
Let $\egyenlet$ be a self-affine set with the  measure separation
property (or in particular with the SSC or SOSC or COSC) 
and let $\mu$ be a self-affine measure on it.
Then for every $\varepsilon>0$ and for every Borel set $B\su K$ with positive
$\mu$-measure
there exists an elementary piece $a(K)$ of $K$ 
of arbitrarily large generation such that
$\m{B\cap a(K)} > (1-\eps) \m{a(K)}$.
\end{lemma}

\bp
Since $\mu(B)>0$, using Fact~\ref{origdef}, 
$B$ can be covered by countably many 
elementary pieces $\phi_{I_i}(K)$ $(i\in\N)$ such that 
$$(1+\eps)\,\mu(B) > \sum_i \m{\phiIiK}.$$
By subdividing the elementary pieces if necessary, we can suppose
that each is of large generation.

If there exists an $i\in\N$ such that
$$(1+\eps)\,\m{B \cap \phiIiK} > \m{\phiIiK}$$
then we can choose $\phiIi$ as $a$.

Otherwise we have
 $(1+\eps)\,\m{B\cap \phiIiK} \le \m{\phiIiK}$ for each $i\in\N$, hence
$$
(1+\eps)\,\mm{B} = 
(1+\eps)\,\m{\bigcup_{i} B\cap \phiIiK} \le 
\sum_i (1+\eps) \m{B\cap \phi_{I_i}(K)} \le
\sum_i \m{\phi_i(K)},
$$ 
contradicting the above inequality.
\qed

\begin{lemma}\label{moho}
Let $\egyenlet$ be a self-affine set with the measure separation
property (or in particular with the SSC or SOSC or COSC) 
and let $\mu$ be a self-affine measure on it.
Then for any Borel set $B\su K$ and $\eps>0$
there exist countably many pairwise almost disjoint elementary pieces $a_i(K)$
such that
$\m{B\cap a_i(K)} > (1-\eps) \m{a_i(K)}$
and 
$\m{B\sm \cuppp_i a_i(K)}=0$.
\end{lemma}

\bp
The elementary pieces $a_i(K)$ will be chosen by greedy algorithm.
In the $n^{th}$ step ($n=0,1,2,\ldots$) 
we choose the largest elementary piece $a_n(K)$ such that
$\m{a_n(K)\cap a_i(K)}=0 \quad (0\le i <n)$ and 
$\m{B\cap a_n(K)}>(1-\eps) \m{a_n(K)}$. 
If there is no such $a_n(K)$ then the procedure terminates.

We claim that $\m{B\sm \cuppp_i a_i(K)}=0$.
Suppose that $\m{B\sm \cuppp_i a_i(K)}>0$.
Then by Lemma~\ref{suruseg} there exists an elementary piece $a(K)$
such that
$$\m{(B\sm \cuppp_i a_i(K)) \cap a(K)} > (1-\eps) \m{a(K)}.$$ 
Then $\m{B\cap a(K)} > (1-\eps) \m{a(K)}$ but $a(K)$ was not chosen
in the procedure.
This could happen only if $a(K)$ intersects a chosen elementary
piece $a_i(K)$ in a set of positive measure. 
But then either $a_i(K)\supset a(K)$ or $a_i(K) \subset a(K)$,
which are both impossible.
\qed

\subsection{Self-affine Sierpi\'nski sponges.}

%
%

\bd
\label{def:K}
By \emph{self-affine Sierpi\'nski sponge} we mean self-affine sets of the
following type.
Let $n,r\in\N$, $m_1,m_2,\ldots,m_n \ge 2$ integers, $M$ be the 
linear transformation given by the diagonal
$n\times n$ matrix 
$$
M=\left(
\begin{array}{ccc}
m_1 &        &   0  \\
    & \ddots &      \\
0   &        & m_n  
\end{array}\right),
$$
and let 
$$
D=\{ d_1,\ldots,d_r\}
\su\{0,1,\ldots,m_1-1\} \times \ldots \times \{0,1,\ldots,m_n-1\}
$$
be given. 
Let $\phi_j(x)=M^{-1}(x+d_j)$ $(j=1,\ldots,r)$ . 
Then the self-affine set $K(M,D)=\egyenlet$ is a Sierpi\'nski sponge.
\ed

We can also define the self-affine Sierpi\'nski sponge as 
$$
K=K(M,D)=\left\{\ \sum_{k=1}^\infty M^{-k} \al_k \ : \ 
\al_1,\al_2,\ldots \in D \ \right\},
$$
or equivalently $K$ is the unique compact set in $\R^n$ 
(in fact, in $[0,1]^n$) such that
$$
M(K)=K+D = \bigcup_{j=1}^r K+d_j;
$$
that is,
$$
K=M^{-1}(K)+M^{-1}(D).
$$
By iterating the last equation we get
\begin{eqnarray*}
  K & = & M^{-k}(K)+M^{-k}(D)+M^{-k+1}(D)+\ldots+M^{-1}(D) \\
    & = & \bigcup_{\al_1,\ldots,\al_k\in D} 
          M^{-k}(K)+M^{-k}\al_k+\ldots+M^{-1}\al_1.
\end{eqnarray*}
Note that the $k$-th generation elementary pieces of $K$ are
the sets of the form
$M^{-k}(K)+M^{-k}(\al_k)+\ldots+M^{-1}(\al_1)$
($\al_1,\ldots,\al_k\in D$) and
the only $0$-th generation elementary part of $K$ is $K$ itself.

\bd\label{standard}
By the \emph{standard} (or sometimes \emph{natural}) 
probability measure on a self-affine sponge
$K=K(M,D)$ we shall mean the self-affine measure on $K$ obtained
by using equal weights $p_j=\frac1r$ ($j=1,\ldots,r$).
\ed

Since the first generation elementary pieces of $K$ are translates
of each other (in fact, so are the $k$-th generation elementary parts),
this is indeed the most natural self-affine measure on $K$.
Using (\ref{origdef}) of Facts~\ref{facts} we get that
$$
\mu(B)=\inf\left\{\sumi \mu(S_i) : B\su\cupi S_i,\ S_i 
\textrm{ is an elementary part of } K\ (i\in\N)\right\}
$$
for every Borel set $B\su K$.

Let $\mut$ be the $\Z^n$-invariant extension of $\mu$ to $\R^n$; 
that is,
for any Borel set $B\su\R^n$ let 
$$
\mut(B)=\sum_{t\in\Z^n}\m{(B+t)\cap K}.
$$

One can check that 
\begin{equation}
\label{mut}
\mut\big(M^l(H) + v\big) = r^l \mu(H) \quad
\textrm{for any } H\su K \textrm{ Borel set, } 
v\in \Z^n,\ l=0,1,2,\ldots.
\end{equation}

\begin{lemma}
\label{approx}
Let $m_1,\ldots,m_n\ge 2$ and $M$ like in Definition~\ref{def:K} 
and let $t\in\R^n$ be 
such that  $\norma{M^k t} >0$ for every $k=0,1,2,\ldots$, where
$\norma{.}$ denotes the distance from $\Z^n$.


Then there exists infinitely many $k\in\N$ such that
$\norma{M^k t} > \frac{1}{2\max(m_1,\ldots,m_n)}.$ 
\end{lemma}

\bp
This lemma immediately follows from the following clear fact:
$$
\norma{u}\le \frac{1}{2\max(m_1,\ldots,m_n)}
\Longrightarrow
\norma{Mu} \ge \min(m_1,\ldots,m_n) \norma{u} \ge 2 \norma{u}.
$$
\qed

\subsection{Invariant extension of measures to larger sets.}

\bl
\label{extendgen}
Suppose that the group $G$ acts on a set $X$, $\iM$ is a $G$-invariant
$\si$-algebra on $X$, $A\in\iM$, $\MA=\{B\in\iM : B\su A\}$
and $\mu$ is a measure on $(A,\MA)$.
 
Then the following two statements are equivalent:

\begin{itemize}
\item[(i)] $\m{g(B)}=\mu(B)$ whenever $g\in G$ and $B,g(B)\in\MA$.
\item[(ii)] There exists a $G$-invariant measure $\mut$ on $(X,\iM)$ such
that $\mut(B)=\mu(B)$ for every $B\in\MA$.
\end{itemize}

\el

\bp
The implication $(ii)\Rightarrow (i)$ is obvious. For proving the other
implication we construct $\mut$ as follows.

If $H$ is a set of the form

\begin{equation}
  \label{partition}
  H=\cupsi B_i, \textrm{ where } g_1,g_2,\ldots\in G 
\textrm{ and } g_1(B_1),g_2(B_2),\ldots\in\MA
\end{equation}
then let 
$$
\mut(H)=\sumi \m{g_i(B_i)}
$$ 
and let $\mut(H)=\infty$ if $H\in\iM$
cannot be written in the above form.


First we check that $\mut$ is well defined; that is, if
we have (\ref{partition}) and
$H=\cupsj C_j$,
$h_1,h_2,\ldots\in G$ and $h_1(C_1), h_2(C_2),\ldots\in\MA$ then
\begin{equation}
  \label{welldef}
  \sumi\m{g_i(B_i)}=\sumj\m{h_j(C_j)}.
\end{equation}
Using that $B_i\su H=\cupsj C_j$ we get that
$g_i(B_i)=g_i(\cupsj B_i\cap C_j) = \cupsj g_i(B_i\cap C_j)$ and so

$$
\sumi\m{g_i(B_i)}=
\sumi\m{\cupsj g_i(B_i\cap C_j)} = \sumi\sumj \m{g_i(B_i\cap C_j)},
$$
and similarly
$$
\sumj\m{h_j(C_j)} = \sumj\sumi \m{h_j(B_i\cap C_j)}.
$$
Thus, using condition (i) for $B=g_i(B_i\cap C_j)$ and 
$g=h_j g_i^{-1}$, we get (\ref{welldef}).

Using the freedom in (\ref{partition}) 
and that whenever $H\in\iM$ can be written in the form
(\ref{partition}) then the same is true for any 
$H\supset H'\in\iM$, 
it is easy to check that 
$\mut$ is a $G$-invariant measure on $(X,\iM)$ such that 
$\mut(B)=\mu(B)$ for every $B\in\MA$.
\qed

We will need only the following special case of this lemma.

\bl
\label{extendspec}
Let $\mu$ be a Borel measure on a Borel set $A\su\Rn$ ($n\in\N$),
$G$ is group of affine transformations of $\Rn$ 
and suppose that
\begin{equation}
  \label{inv}
  \m{g(B)}=\mu(B) \textrm{ whenever } b\in G,\ B, g(B)\su A 
  \textrm{ and } B \textrm{ is a Borel set}. 
\end{equation}

Then there exists a $G$-invariant Borel measure $\mut$ on
$\Rn$ such that $\mut(B)=\mu(B)$ for any $B\su A$ Borel set. 
\ep
\el

\br
The extension we get in the above proof do not always give the
measure we expect -- it may be infinity for too many sets.
For example, if $A\su\R$ is a Borel set of first category with
positive Lebesgue measure, $G$ is the group of translations
and $\mu$ is the restriction of the 
Lebesgue measure to $A$ then the Lebesgue measure itself would be
the natural translation invariant extension of $\mu$,
however the extension $\mut$ as defined in the proof is clearly
infinity for every Borel set of second category.
\er

\begin{ddd}\label{isominvdef}\upshape
Let $\mu$ be a Borel measure on a compact set $K$. We say that $\mu$ is isometry invariant
if given any isometry $g$ and a Borel set $B\subset K$ such that $g(B)\subset K$,
then $\mu(B)=\m{g(B)}$.
\end{ddd}

This definition makes sense since (by Lemma~\ref{extendspec}) exactly
the isometry invariant measures on $K$ can be extended to be isometry invariant measures
on $\R^n$ in the usual sense.

As an illustration of Lemma~\ref{extendspec} we mention the following special case with a
peculiar consequence.

\bl
\label{specspec}
Let $A\su\Rn$ ($n\in\N$) be a Borel set such that $A\cap(A+t)$
is at most countable for any $t\in\Rn$. Then any continuous 
Borel measure $\mu$ on $A$ (continuous here means that the measure
of any singleton is zero) can be extended to a translation invariant
Borel measure on $\Rn$. 
\ep
\el

Note that although the condition that $A\cap(A+t)$
is at most countable for any $t\in\Rn$ seems to imply that $A$ is very
small, such a set can be still fairly large. 
For example there exists a compact set $C\su\R$ 
with Hausdorff dimension $1$ such that
$C\cap (C+t)$ contains at most one point for any $t\in\R$ \cite{Ke}.
Combining this with Lemma~\ref{specspec} we get the following.

\bcoro
\label{dimone}
There exists a compact set $C\su\R$ 
with Hausdorff dimension $1$ such that
any continuous 
Borel measure $\mu$ on $C$ can be extended to a translation invariant
Borel measure on $\R$.
\ep
\ecoro

\subsection{Some more lemmas.}

The following simple lemmas might
be known but for completeness
(and because it is easier to prove them than to find them) we present
their proof.

Recall that the support of a measure is the smallest closed set
with measure zero complement.

\begin{lemma}
\label{continuity}
Let $\mu$ be a finite Borel measure on $\R^n$ 
with compact support $K$.
Then for every $\eps>0$ there exists a $\de>0$ such that
$$
|u|\ge\eps \Longrightarrow \m{K\cap (K+u)} \le (1-\de)\mu(K).
$$
\end{lemma}


\bp
We prove by contradiction. Assume that there exists an 
$\eps>0$ and a sequence $u_1,u_2,\ldots\in\R^n$ 
such that $|u_n|\ge \eps$ (for every $n\in\N$)
and $\m{K\cap(K+u)}\to\mu(K)>0$ $(n\to\infty)$.
By omitting some (at most finitely many) zero terms
we can guarantee that
every $u_n$ is in the compact annulus
$\{x : \eps \le |x| \le \diam(K)\}$
(where $\diam$ denotes the diameter),
so by taking a subsequence we can suppose
that $(u_n)$ converges, say to $u$.
Since $K\cap (K+u)$ is a proper compact subset of $K$
(since $K$ is compact and $u\neq 0$, $K+u\supset K$
is impossible) and $K$ is the support of $\mu$,
we must have $\mu(K)>\m{K\cap (K+u)}=\mu(K+u)$.

It is well known (see e.g. \cite{Ru}, 2.18. Theorem)
that any finite Borel measure is outer regular
in the sense that the measure of any Borel set
is the infimum of the measures of the open sets that contain
the Borel set.
Thus $\mu(K+u)<\mu(K)$ implies that
there exists an open set 
$G \supset K+u$ such that 
$\mu(G)<\mu(K)$.
Then whenever $|u_n-u|$ is less than the (positive)
distance between $K$ and the complement of $G$, $G$ 
contains $K+u_n$ 
and so $\mu(K)>\mu(G)\ge\mu(K+u_n)$.
This is a contradiction since $u_n\to u$
and $\mu(K+u_n)=\m{K\cap(K+u_n)}\to\mu(K)$.
\qed

\begin{lemma}\label{fff}
Let $K\su\Rd$ be compact and $\mu$ be a probability Borel measure
on $K$ such that any nonempty relative open subset of $K$ has positive
$\mu$ measure. 
Then if the sequence $(g_n)$ of affine maps
converges to an affine map $g$ and $\m{g_n(K)\cap K} \to 1$
then $\m{g(K)\cap K}=1$. Moreover, $K \subset g(K)$.
\end{lemma}

\bp
Suppose that $\m{g(K)\cap K}=q<1$. 
Let $g(K)_\eps$ denote the $\eps$-neighborhood of $g(K)$.
Since $\bigcap_{n=1}^\infty (g(K)_{1/n} \cap K) = g(K) \cap K$
and $\mu$ is a finite measure we have 
$\m{g(K)_{1/n} \cap K} \to \m{g(K)\cap K}=q$. 
Thus there exists an $\eps>0$ for which 
$\m{g(K)_\eps \cap K} \le \frac{1+q}{2} <1$. 
Since $g_n$ converges uniformly on $K$, for $n$ large enough we have
$g_n(K) \subset g(K)_\eps$ and so
$\m{g_n(K)\cap K} \le \frac{1+q}{2}$,
contradicting $\m{g_n(K)\cap K} \to 1$.
Therefore we proved that $\m{g(K)\cap K}=1$.

Then $K\sm g(K)$ is relative open in $K$
and has $\mu$ measure zero, so it 
must be empty, therefore $K\su g(K)$.
\qed

\section{Self-affine sets with the strong separation condition}
\label{affine}

\begin{propo}\label{taraza}
For any self-affine set $K\su\Rd$ with the strong separation condition
there exists an open neighborhood $U\su \Ak$ of the identity map 
such that 
for any $g\in U$, $$g(K)\supset K \Longleftrightarrow g=\ii.$$
\end{propo}

\bp
Let $n$ denote the dimension of the affine span of $K$.

We shall prove that there exists 
a small open neighborhood $V\su \Ak$ of the identity map
such that for any $g\in V$ we have $g(K)\subset K \Longleftrightarrow g=\ii$.  
This would be enough since then for any $g\in V$
we get $K \subset g^{-1}(K) \Longleftrightarrow g =\ii$, therefore
$U=V^{-1}=\{g^{-1}: g\in V\}$ has all the required properties.

Similarly as in the proof of Lemma~\ref{jancski}, choose 
$n+1$ elementary pieces $\phi_{I_1}(K), \ldots, \phi_{I_{n+1}}(K)$ of $K$
so that if we pick one point from the convex hull of each of them
then we get a maximal collection of affine independent points
in the affine span of $K$.

Let $d=\min_{1\le i\le n+1} \dist(\phiIiK, K\setminus\phiIiK)$, then $d>0$.
Let $V$ be a so small neighborhood of the identity map that
$\dist(x,g(x)) < d$ for any $g\in V$ and $x\in K$.

Let $g\in V$ and $g(K)\subset K$. 
Then, by the definition of $d$ and $V$ we have 
$g(\phiIiK)\subset \phiIiK$ for every $1\le i\le n+1$.
Then the convex hulls of these elementary pieces are also mapped
into themselves. Since each of these convex hulls is homeomorphic to 
a ball, by Brouwer's fixed point theorem we get a fixed point of $g$
in each of these elementary pieces. So we obtained $n+1$
fixed points of $g$ such that their affine span is exactly the
affine span of $K$. Since $g$ is an affine map,
the set of its fixed points form an affine subspace, thus the set
of fixed points of $g$ contains the affine span of $K$. 
Since $g\in\Ak$, $g$ is defined exactly on 
the affine span of $K$, therefore $g$ must be the identity map.
\qed

\begin{theorem}\label{kovesarki}
Let $\egy$ be a self-affine set satisfying the strong separation condition
and let $\mu$ be a self-affine measure on $K$. 
Then there exists a $c<1$ and an open neighborhood $U\su\Ak$ of the
identity map such that 
$g\in U\sm\id \Longrightarrow \m{K\cap g(K)} < c$.
\end{theorem}

\bp
Using Proposition~\ref{taraza} we can choose a small open 
neighborhood $U\su\Ak$ of the
identity map such that even in the closure of $U$ the only affine map $g$
for which $g(K)$ contains $K$ is the identity map and so that 
\begin{equation}\label{e:notfar}
\dist(x,g(x)) < 1 \textrm{ for any } g\in U \textrm{ and } x\in K.
\end{equation}
Since $\Ak$ is locally compact, we may also assume that the closure
of $U$ is compact.

We claim that we can choose an
even smaller open neighborhood $V\su U$ of the identity map
such that $\phi_i^{-1}\circ V\circ \phi_i \subset U$ 
for $i=1,\ldots,r$ and that 
$g(\phi_i(K)) \cap \phi_j(K)= \emptyset$
for any $i\neq j$ and $g\in V$.
Indeed, the first property can be satisfied since $\Ak$ is 
a topological group and those $g$'s for which the second property 
do not hold are far from the identity map.

Now we claim that there exists a $c<1$ such that
$g\in \ubolv \Longrightarrow \m{g(K)\cap K}<c$.
Suppose that there exists a sequence $(g_n)\su\ubolv$ such that
$\m{K\cap g_n(K)}\to 1$. 
Since $\ubolv$ is compact there exists a subsequence $g_{n_i}$ such that
 $g_{n_i}\to h\in \ubolv$.
By Lemma~\ref{fff} this implies that $h(K)\supset K$
but in $\ubolv$ there is no such affine map $h$.

We prove that this $U$ and this $c$ have the required properties;
that is, 
$g\in U\sm\id\Longrightarrow\m{K\cap g(K)}< c$.

If $g\in\ubolv$ then we are already done, so suppose that $g\in V\sm\id$. 
Let $F$ denote the set of fixed points of $g$.

The heuristics of the remaining part of the proof is the following.
The affine map $g$ moves $K$ too slightly. We zoom in on small
elementary pieces $a(K)$ of $K$ so that each $g(a(K))$ intersects only $a(K)$
in $K$, but $g$ moves $a(K)$ far enough (compared to its size).
Technically this second requirement means that 
$a^{-1}\circ g \circ a\in U\sm V$, so we can use the $g\in U\sm V$ case
for the elementary piece $a(K)$.
We find such an elementary piece around each point of $K$ that is
not a fixed point of $g$, and so we get
a partition of $K\sm F$
into elementary pieces with the above property.
Finally, by adding up the estimates for these elementary pieces we derive 
$\m{g(K)\cap K}<c$.

\bc\label{c1}
For any $x\in K\sm F$ there exists a largest
elementary piece $\phi_{I_x}(K)$ of $K$ that contains $x$ and for which
$\phi_{I_x}^{-1} \circ g \circ \phi_{I_x} \in U\sm V$.
\ec
\bp
Let $(i_1, i_2, \ldots)$ be the sequence of indices for which 
$$\{x\}=\bigcap_{n=1}^\infty (\phi_{i_1} \circ \phi_{i_2} 
\circ \ldots \circ \phi_{i_n})(K),
$$
and let $I_n=(i_1, \ldots, i_n)$. 
Since $g\in V$, we have $\phi_{i_1}^{-1} \circ g \circ \phi_{i_1} \in U$
by the definition of $V$.
If for some $n$ we have 
$\phi_{I_n}^{-1} \circ g \circ \phi_{I_n} \in V$
then by the definition of $V$ we have 
$$
\phi_{I_{n+1}}^{-1} \circ g \circ \phi_{I_{n+1}} =
\phi_{i_{n+1}}^{-1} \circ \phi_{I_n}^{-1} \circ g 
\circ \phi_{I_n} \circ \phi_{i_{n+1}} \in U.
$$
Therefore it is enough to find an $n$ such that 
$\phi_{I_n}^{-1} \circ g \circ \phi_{I_n} \not\in V$
since then taking the smallest such $n$, $I_x=I_n$ has the 
desired property.
Letting $y_n=\phi^{-1}_{I_n}(x)$ we have $y_n\in K$ 
(since $\{x\}=\bigcap_{n=1}^\infty \phi_{I_n}(K)$) and 
$(\phi^{-1}_{I_n} \circ g \circ \phi_{I_n})(y_n)=\phi_{I_n}^{-1}(g(x))$.
Since $x$ is not a fixed point of $g$, for $n$ large enough we have
$$
\dist\big(g(x), \phi_{I_n}(K)\big) > \frac{\dist(g(x),x)}{2} \= t >0.
$$
Since each $\phi_i$ is a contractive affine map, there exists an $\al_i<1$
such that $\dist(\phi_i(a),\phi_i(b))\le \al_i \cdot \dist(a,b)$ 
for any $a,b$.
Then, using the multi-index notation
$\al_{I_n}=\al_{i_1}\cdot\ldots\cdot\al_{i_n}$, we clearly have 
$\dist(\phi_{I_n}(a),\phi_{I_n}(b))\le \al_{I_n}\cdot\dist(a,b)$ 
for any $a,b$.
Then
$\dist\big(\phi^{-1}_{I_n}(g(x)), K\big) > t/\alfa_{I_n}$, hence 
$\dist\big((\phi^{-1}_{I_n}\circ g \circ \phi_{I_n})(y_n), K\big) > 
t/\alfa_{I_n}$, 
which is bigger than $1$ if $n$ is large enough.
Thus for $n$ large enough, $\phi^{-1}_{I_n} \circ g \circ \phi_{I_n}$
is not in $V$, since it is not even in $U$ by (\ref{e:notfar}).
\qed

\bc \label{c2}
For any $x\in K\sm F$ we have
$g(\phi_{I_x}(K)) \cap K \subset \phi_{I_x}(K)$, where
$I_x=I_n=(i_1,\ldots,i_n)$ is the multi-index we got in Claim~\ref{c1}.
\ec
\bp
Let $k\in\{0,1,\ldots,n-1\}$ be arbitrary and let $I_k=(i_1,\ldots,i_k)$.
Then $\phii_{I_k} \circ g \circ \phi_{I_k} \in V$, hence 
for any $l \noteq i_{k+1}$ we have
$(\phii_{I_k} \circ g \circ \phi_{I_k} \circ \phi_{i_{k+1}})(K) \cap \phi_l(K) 
= \ures$, which is the same as
$(g \circ \phi_{I_{k+1}})(K) \cap (\phi_{I_k} \circ \phi_l)(K) 
= \ures$ ($l\noteq i_{k+1}$). 
Since $(g \circ \phi_{I_n})(K)\su (g \circ \phi_{I_{k+1}})(K)$, this
implies that
$$
(g \circ \phi_{I_n})(K) \cap (\phi_{I_k} \circ \phi_l)(K) = \ures \quad
(k\in\{0,1,\ldots,n-1\},\ l\noteq i_{k+1}).
$$
Since $K\sm \phi_{I_n}(K)=\cup_{k=0}^{n-1}\cup_{l\noteq i_{k+1}} 
(\phi_{I_k} \circ \phi_l)(K)$,
this implies that
$g(\phi_{I_n}(K)) \cap K \subset \phi_{I_n}(K)$.
\qed

The elementary pieces $\{\phi_{I_x}(K): x\in K\sm F\}$
clearly cover $K\sm F$.
Since for any $x\neq y$ we have 
$\phi_{I_x}(K) \cap \phi_{I_y}(K) = \ures$ or 
$\phi_{I_x}(K) \subset \phi_{I_y}(K)$ or $\phi_{I_x}(K) \supset
\phi_{I_y}(K)$,
one can choose a
\begin{equation}\label{e:1.5}
K\sm F \subset \bigcupp_{i=1}^{\infty} \phi_{J_i}(K)
\end{equation}
countable disjoint subcover.
By Claim~\ref{c2} we have 
\begin{equation}\lab{e:parammm}
g(\phi_{J_i}(K)) \cap K\subset \phi_{J_i}(K).
\end{equation}

Since $g$ is not the identity map (of the affine span of $K$) 
and $F$ is the set of fixed points of the affine map $g$, the dimension
of the affine subspace $F$ is smaller than the dimension of the affine
span of $K$, and so we cannot have $g(F)\supset K$. 
By Lemma~\ref{jancski} this implies that $\m{g(F)\cap K}=0$.
Using this last equation, (\ref{e:1.5}), (\ref{e:parammm}), 
and finally the definition of a self-affine measure we get that
\begin{multline}
\m{g(K)\cap K} \le 
\m{g(F) \cap K} + \m{g(K\sm F) \cap K} = 
\m{g(K\sm F) \cap K} \\
\le \M{g\Big(\bigcupp_{i=1}^{\infty} \phi_{J_i}(K)\Big) \cap K}= 
\sum_{i=1}^\infty \M{g\big(\phi_{J_i}(K)\big) \cap K} \\
=\sum_{i=1}^\infty \M{g\big(\phi_{J_i}(K)\big) \cap \phi_{J_i}(K)} = 
\sum_{i=1}^\infty \M{\phi_{J_i}\big((\phii_{J_i} \circ g \circ \phi_{J_i})(K)\cap K\big)} \\
=\sum_{i=1}^\infty p_{J_i}\,\M{(\phii_{J_i} \circ g \circ \phi_{J_i})(K)\cap K}.
\nonumber
\end{multline}

Since $\phii_{J_i} \circ g \circ \phi_{J_i} \in U\sm V$,  
the measures in the last expression are less than $c$.
Thus  $\m{g(K) \cap K} < \sum p_{J_i} \cdot c = \sum \m{\phi_{J_i}(K)}\cdot
c = \m{\bigcup^{*}_i \phi_{J_i}(K)}\cdot c =c$, which completes
the proof.
\qed

\begin{theorem}\label{selfaffine}
Let $K\su\Rd$ be a self-affine set with the strong separation condition
and let $\mu$ be a self-affine measure on $K$. 
Then there exists a constant $c<1$ such that for any isometry $g$ 
we have $\m{K\cap g(K)}<c$ unless $g(K)=K$.
\end{theorem}

\bp
Suppose that $g_n\in\Ik$ 
(that is, $g_n$ is an isometry of the affine span of $K$)
such that $g_n(K)\neq K$ ($n\in\N$)
and $\m{K\cap g_n(K)}\to 1$. 
We can clearly assume that $K\cap g_n(K)\neq\emptyset$ for each $n$ and
so the whole sequence $(g_n)$ is in a compact subset of $\Ik$.
Thus, after choosing a subsequence if necessary, 
we can also assume that $g_n$ converges to an $h\in\Ik$.
By Lemma~\ref{fff} we must have $K\su h(K)$.
It is well known and not hard to prove that no compact set in $\Rd$ can
have an isometric proper subset, so $K\su h(K)$ implies that $h(K)=K$.

Applying Theorem~\ref{kovesarki} we get a $c<1$ and an open 
neighborhood  $U\su\Ak$ of the identity such that
$g\in U\sm\id \Longrightarrow \m{K\cap g(K)}<c$.

Since $g_n\to h$ we get $g_n\circ h^{-1}\to \ii$.
Let $n$ be large enough to have $g_n\circ h^{-1}\in U$ 
and $\m{K\cap g_n(K)}>c$.
Since $g_n(K)\neq K$ but $h(K)=K$ we cannot have $g_n=h$ and so 
$g_n\circ h^{-1}\in U\sm\id$. Then, by the previous paragraph, we get
$\m{K\cap g_n(K)}<c$, contradicting $\m{K\cap g_n(K)}>c$.
\qed

\section{Self-similar sets with the strong separation property}
\label{similar}

Our first goal in this section is to prove the following theorem.

\begin{theorem}\label{andras}
Let $\egy$ be a self-similar set satisfying the strong separation condition 
and $\mu$ be a self-similar measure on it. 
There exists $c<1$ such that for every similitude
$g$ either $\m{g(K)\cap K}<c$ or $K\subset g(K)$.
\end{theorem}

Now, for the sake of transparency we outline the proof. At first we need a new notation.

\def\Skcs{S^*_K}
From $S_K$ we excluded those similarity maps which map everything to a single point. So let $\Skcs$ be the metric space of all degenerate and all non-degenerate similarity maps in the affine span $A_K$ of $K$; that is,
\begin{equation}\label{eq:skcs}
\Skcs=S_K \cup \{f\ |\ f:A_K\to \{y\},\, y\in A_K\}.
\end{equation}

First we show that there exists a compact set $\G\su\Skcs$ of similarity maps such
that for every $g'\in \Skcs$ there exists $g\in \G$ such that $g'(K)\cap
K=g(K)\cap K$. 
Then it is easy to see that it suffices to prove the theorem for $g\in\G$.
(It is easy to see that no such compact set $G$ in $\Sk$ exists.)

Let $\mu_H$ be a constant multiple of Hausdorff measure of appropriate
dimension so that $\mu_H(K)=1$.
The restriction of this measure to $K$ is a self-similar measure.
Let us consider those $h\in \G$ for which $K\subset h(K)$ holds.
Using Hausdorff measures and Theorem~\ref{kovesarki} we prove that there are
only finitely many such $h$, and also that the theorem holds in small
neighbourhoods of each such $h$ for the measure $\mu_H$. The maximum of the
corresponding finitely many values $c$ is still strictly smaller than $1$. Let
us now cut these small neighbourhoods out of $\G$. Using upper semicontinuity
of our measure (Lemma~\ref{fff}) we produce a $c<1$ such that for the
remaining similarity maps $g$ we have $\mu_H\big(g(K)\cap K\big)<c$. 
Then clearly the same holds for all elements of $\G$, possibly with a larger
$c<1$, finishing the proof for the measure $\mu_H$.

Applying the theorem for $\mu_H$, and also in a small open neighbourhood $U$
of the identity for every self-similar measure $\mu$, we show that if $h\in\G$,
$K\subset h(K)$, and $g$ is in a small neighbourhood of $h$ then $\m{g(K)\cap
K}<c$. Then the same argument as above (using upper semicontinuity) yields the
theorem, possibly with a larger constant again.

\begin{propo}\label{korlatos}
Let $\egy$ be a self-similar set satisfying the strong separation condition. Then
there exists a compact set $\G \subset \Skcs$ such that for every 
similarity map
$g'\in\Skcs$ there is a $g\in \G$ for which $g'(K)\cap K = g(K) \cap K$ holds.
\end{propo}

\bp
Let $D$ denote the diameter of $K$, let 
$\delta=\min_{1\le i < j \le r} \dist(\phi_i(K), \phi_j(K))$ and let
$$
\G=\{g\in\Skcs : g(K)\cap K\neq\emptyset, \textrm{ the similarity ratio of }
g \textrm{ is at most } D/\de\}\cup\{g_0\},
$$
where $g_0\in\Skcs$ is an arbitrary fixed similarity map such that  
$g(K)\cap K=\emptyset$.
It is easy to check that $\G\su\Skcs$ is compact.

Let $g'\in\Skcs$. If $g'\in\G$ or $g'(K)\cap K=\emptyset$ then we can choose
$g=g'$ or $g=g_0$, respectively. So we can suppose that 
$g'(K)\cap K\neq\emptyset$ and the similarity ratio of 
$g'$ is greater than $D/\delta$. Then the
minimal distance between the first generation elementary pieces $g'(\phi_j(K))$
of $g'(K)$ is larger than $D$. 
So there exists $\phi_i$ such that $g'(K)\cap K = g'(\phi_i(K))\cap K$.
Therefore $g'$ can be replaced by $g'\circ \phi_i$, which has similarity
ratio $\al_i$ times smaller than the similarity ratio of $g'$, where $a_i$
denotes the similarity ratio of $\phi_i$. 
Since $\max(\al_1,\ldots,\al_r)<1$, this way in finitely many steps we
get a $g$ with similarity ratio at most $D/\de$ such that 
$g(K)\cap K = g'(K)\cap K\neq\emptyset$, which completes the proof.
\qed

\begin{propo}\label{haus}
Let $\egy$ be a self-similar set satisfying the strong separation condition. 
\begin{enumerate}[(i)]
\item\lab{egy} 
Then $\{g\in\Sk : g(K)\supset K\}$ is discrete in $\Sk$, hence countable,
and also closed in $\Sk$.
\item\lab{ketto} 
Let $\mu_H$ be a constant multiple of Hausdorff measure of appropriate
dimension so that $\mu_H(K)=1$. There exists $c<1$ such that for every
similitude $g$ either
$\mu_H\big(g(K)\cap K\big)<c$ or $K \subset g(K)$.
\end{enumerate}
\end{propo}

\bp
By Lemma~\ref{fff} $\{g\in\Sk : g(K)\supset K\}$ is closed. 
Since every discrete subset of a subspace of $\R^{d^2+d}$ is
countable, in order to prove (\ref{egy}) it is enough to prove that
$\{g\in\Sk : g(K)\supset K\}$ is discrete.

Let $\eps$ be a positive number to be chosen later, 
and $h$ be a similitude for which $K\subset
h(K)$. Denote by $K_\delta$ the $\delta$-neighbourhood of $K$. As 
$\muH\big(h(K)\big)$ is finite, there is a small $\delta>0$ such that
$\muH\big(K_\delta \cap (h(K)\sm K)\big)<\eps$.  Applying
Theorem~\ref{kovesarki} to $K$ and $\mu_H$ we obtain an open neighbourhood $U\su\Ak$
and a constant $c_H$. There exists an open neighbourhood $W_\varepsilon\subset
\Sk$ of the identity such that
\begin{itemize}
\item[(a)] $W_\varepsilon = W_\varepsilon^{-1} \subset U$,
\item [(b)] $\dist(g(x),x)<\delta$ for every $x\in K$,
\item[(c)] $\muH\big(g(B)\big) \le
(1+\eps)\muH(B)$ for every $g\in W_\varepsilon$ and Borel set $B$,
\end{itemize}
where for (c) we use that a similitude of ratio $\al$ multiplies 
the $s$-dimensional Hausdorff measure by $\al^s$.

Let $g\in W_\varepsilon h$ and $g\neq h$. 
Clearly $W_\varepsilon h$ is an open neighbourhood of $h$ and
$g\circ
h^{-1}$, $h\circ g^{-1} \in W_\varepsilon\minid$, and $(h\circ
g^{-1})(K)\subset K_\delta$.  Hence
\begin{multline}\label{long}
\muH\big(K\cap g(K)\big) \le 
(1+\eps) \muH\big((h\circ g^{-1})(K \cap g(K))\big) = \\
= (1+\eps) \muH\big((h\circ g^{-1})(K) \cap h(K)\big)= \\ 
=(1+\eps) \muH\big( (h\circ g^{-1})(K) \cap K\big) +(1+\eps)\muH\big((h\circ
g^{-1})(K) \cap (h(K)\sm K)\big) \le \\ 
\le (1+\eps)c_H + (1+\eps)\muH\big(K_\delta \cap (h(K)\sm K)\big) 
\le (1+\eps) c_H + (1+\eps) \eps.
\end{multline}
The last expression is clearly smaller than $1$ if $\eps$ is small enough, so
let us fix such an $\eps$. Therefore if $g\in W_\varepsilon h$ and $g\neq h$
then $g(K)\not\supset K$, which shows that $\{g\in\Sk : g(K)\supset K\}$ is
discrete finishing the proof of
(\ref{egy}).

In order to prove (\ref{ketto}) suppose towards a contradiction 
that $\sup\,\{\muH(g(K)\cap K) :
g\in\Skcs, \,g(K)\not\supset K\}=1$. Then we also have
$\sup\,\{\muH(g(K)\cap K) : g\in\G,\,g(K)\not\supset K\}=1$. Let
$(g_n)$ be a
convergent sequence in $\G$ so that $g_n(K)\not\supset K$, $\muH\big(g_n(K)\cap
K\big)\to 1$, $g_n \to h$. Lemma~\ref{fff} yields $h(K)\supset K$, hence $g_n \neq
h$. If $n$ is large enough then $g_n\in W_\varepsilon h$ and, by (\ref{long}),

$\muH\big(K\cap g_n(K)\big)\le
(1+\eps)c_H+(1+\eps)\eps$, contradicting $\muH\big(g_n(K)\cap K\big)\to 1$.
\qed

%


\proc{Proof of Theorem~\ref{andras}.}
By Proposition~\ref{korlatos} we can assume $g\in\G$.
Let $c_H$ be the constant yielded by Proposition~\ref{haus}~(\ref{ketto}).
Fix $h\in \G$ with $h(K)\supset K$. There are only finitely many such $h$ by
Proposition~\ref{haus}~(\ref{egy}) and the compactness of $\G$.

Let us now apply Lemma~\ref{suruseg} to the self-similar set $h(K)$,
$\mu_H$, $0<\eps\le 1-c_H$ and $B=K\subset h(K)$.
We obtain $\phi_I$ such that 
\[
\mh{K\cap h(\phi_I(K))} \ge (1-\eps)\,\mh{h(\phi_I(K))}.
\]
Hence Proposition~\ref{haus}~(\ref{ketto}) applied to 
the self-similar set $h(\phi_I(K))$ and the similitude $(h\circ \phi_I)^{-1}$ 
gives $K\supset h(\phi_I(K))$.

Since $h(\phi_I(K))$ is open in $h(K)$, it is also open in $K$
and so it can be written as a union of elementary pieces of $K$. 
Since $h(\phi_I(K))$ is compact this implies
that $h(\phi_I(K))$ is a finite union of elementary
pieces of $K$. Let $\phi_J(K)$ be one of these elementary pieces. 
So $\phi_J(K) \subset
h(\phi_I(K)) \subset K \subset h(K)$. As $\phi_J(K)$ is open in $K$,
it is also open in $h(\phi_I(K))$, hence also in $h(K)$. Therefore $\dist(
\phi_J(K), h(K)\sm \phi_J(K)) > 0$, 
and so for every $g$ that is close enough to
$h$ we have
\[
\big(g\circ h^{-1}\big)\big(h(K)\sm \phi_J(K)\big) \cap \phi_J(K) = \ures.
\]
Thus, as $\phi_J(K)\su h(K)$, for every such $g$ we have
$$
g(K)\cap \phi_J(K) = \big(g\circ h^{-1}\big)\big(h(K)\big) \cap \phi_J(K) 
= \big(g\circ h^{-1}\big)\big(\phi_J(K)\big) \cap \phi_J(K).  
$$
On the other hand, Theorem~\ref{kovesarki} yields that 
there exists a $c<1$ such that if $g$ is close enough to $h$ and  
$g\neq h$ then
\[
\m{(g\circ h^{-1})(\phi_J(K)) \cap
\phi_J(K)} < c \cdot \mu\big(\phi_J(K)\big)= c\cdot p_J.
\]
Therefore
$\m{g(K)\cap \phi_J(K)}=\m{(g\circ h^{-1})(\phi_J(K)) \cap \phi_J(K)} < c\cdot
p_J$ and
\begin{multline}
\m{g(K)\cap K} = \m{g(K)\cap \phi_J(K)} + \m{g(K) \cap (K\sm \phi_J(K))}\\
< c\cdot p_J + 1-p_J=1-(1-c)p_J.
\end{multline}
As we only considered finitely many $h$'s,
there exists $c'<1$ such that if $g$ is close to one of these $h$'s, but
distinct from it, then $\m{g(K)\cap K}<c'$. This, together with
Lemma~\ref{fff} provides a $c''<1$ such that for every $g\in \G$ either
$\m{g(K)\cap K} < c''$ or $g(K)\supset K$. (Just like at the end of the proof of
Proposition~\ref{haus}.) 
Finally, by Proposition~\ref{korlatos} this also holds outside $\G$.
\qed

We will apply this theorem to elementary pieces of $K$ instead of $K$
itself. It is easy to see that the same $c$ works for every elementary
piece; that is, we have the following corollary of Theorem~\ref{andras}.

\begin{cccc}\label{cor}
Let $\egy$ be a self-similar set satisfying the strong separation condition 
and $\mu$ be a self-similar measure on it. 
There exists $c<1$ such that for every similitude
$g$ and every elementary piece 
$a(K)$ of $K$ either
$\mu\big(g(K) \cap a(K)\big)< c\cdot \mu\big(a(K)\big)$ 
or $a(K) \subset g(K)$.
\ep
\end{cccc}

Now we are ready to prove the second main result of this section.

\begin{theorem}\label{pos}
Let $\egy$ be a self-similar set satisfying the strong separation condition, $\mu$
be a self-similar measure on it, and $g$ be a similitude. Then $\m{g(K)\cap K}
>0$ if and only if the interior (in $K$) of $g(K)\cap K$ is nonempty.
Moreover, $\m{\int_K (g(K)\cap K)}=\m{g(K)\cap K}$.
\end{theorem}

\bp If the interior (in $K$) of $g(K)\cap K$ is nonempty then clearly it is of
positive measure, since the measure of every elementary piece is positive.

Let $c$ be the constant given by Corollary~\ref{cor}, and let $g$ be a
similitude such that $\m{g(K)\cap K} >0$.  Applying Lemma~\ref{moho} 
for $B=g(K)\cap K$ and 
$\eps=1-c$ we obtain countably many disjoint elementary pieces $a_i(K)$ of $K$
such that
\begin{equation}\label{e:big}
\m{g(K)\cap a_i(K)}= \m{(g(K)\cap K) \cap a_i(K)} > c \cdot
\m{a_i(K)}
\end{equation}
and $\big(g(K)\cap K\big)\sm \bigcup^*_i a_i(K)$
is of $\mu$-measure zero. 
By Corollary~\ref{cor}, (\ref{e:big}) implies that $a_i(K) \subset g(K)$. 
Since $a_i(K)$ is open in $K$, it is open in $g(K)\cap K$, so $\bigcup^*_i
a_i(K) \subset \int_K (g(K)\cap K)$. Hence 
\mult{
\m{g(K)\cap K} = \m{g(K)\cap K \cap \bigcupp_i a_i(K)} + \m{(g(K)\cap K)\sm
\bigcupp_i a_i(K)}\\ = \m{\bigcupp_i a_i(K)} 
\le \m{\int_K(g(K)\cap K)},
} 
proving the theorem.
\qed

As an immediate consequence we get the following.

\begin{cccc}\label{mindegy}
Let $K\su\Rd$ be a self-similar set satisfying the strong separation condition,
and let $\mu_1$ and $\mu_2$ be self-similar measure on $K$.
Then for any similitude $g$ of $\R^d$,
$$
\mu_1\big(g(K)\cap K\big)>0 \Longleftrightarrow \mu_2\big(g(K)\cap K\big)>0.
$$
\end{cccc}

We also get the following fairly easily.

\begin{cccc}\label{countable}
Let $K\su\Rd$ be a self-similar set satisfying the strong separation condition,
let $\affk$ be the affine span of $K$
and let $\mu$ be a self-similar measure on $K$.
Then the set of those similitudes $g:\affk\to\Rd$ 
for which $\m{g(K)\cap K}>0$ is countably infinite.
\end{cccc}

\bp
It is clear that there exist infinitely many similitudes $g$ such
that $\m{g(K)\cap K}>0$ since the elementary pieces of $K$ are similar
to $K$ and have positive $\mu$ measure.

By Lemma~\ref{jancski}, $\m{g(K)\cap K}>0$ implies that $g\in\Sk$ and,
by Theorem~\ref{pos}, that $g(K)$ contains an elementary piece of $K$. 
Therefore it is
enough to show that for each fixed elementary piece $a(K)$ of $K$
there are only countably many $g\in\Sk$ such that $g(K)\supset a(K)$,
which is the same as $(a^{-1}\circ g)(K)\supset K$. 
By 
the first part of Proposition~\ref{haus} there are only countably many such 
$a^{-1}\circ g\in\Sk$,
so there are only countably many such $g\in\Sk$.
\qed


From the first part of Proposition~\ref{haus} we get more results about those
similarity maps that map a self-similar set into itself. These results will be
used in the next section and they are also related to a theorem and a question
of Feng and Wang \cite{FW} as it will be explained before Corollary~\ref{log}.

\begin{lemma}\label{vegessok}
Let $\egy$ be a self-similar set with strong separation condition. There exists only finitely many similitudes $g$ for which
$g(K) \subset K$ holds and $g(K)$ intersects at least two first generation elementary pieces of $K$.
%
\end{lemma}
\Biz
The similarity ratios of these similitudes $g$ are strictly separated from zero. Thus the similarity ratio of their inverses have some finite upper bound,
and also $K\subset g^{-1}(K)$ holds. The set of similitudes with the latter property form a discrete and closed set according to the first part of
Proposition~\ref{haus}.

%
Those $h\in\Skcs$ similarity maps (cf. \eqref{eq:skcs}) whose similarity ratio is under some fixed bound and for which $h(K)\cap K \neq \emptyset$ holds form a compact set in $\Skcs$
(see proof of Proposition~\ref{korlatos}). Since a discrete and closed subspace of a compact set is finite, the proof is finished.
\qed

\begin{theorem}\label{lambdatul}
Let $\egy$ be a self-similar set with strong separation condition and let $\l$ be a
similitude for which $\l(K)\subset K$. There exist
an integer $k \ge 1$ and multi-indices $I, J$ such that $\l^k \circ \phi_I = \phi_J$.
%
\end{theorem}
\Biz
For every integer $k\ge 1$ there exists a smallest elementary piece $\phiIK$ which contains $\l^k(K)$.
For this multi-index $I$, $(\phii_I \circ \l^k)(K)$ is a subset of $K$ and intersects at least two first generation elementary pieces of $K$.
There are only finitely many similitudes with this property according to Lemma~\ref{vegessok}, hence there exist $k<k'$, $I$, $I'$ such that 
$\phii_I \circ \l^k = \phii_{I'} \circ \l^{k'}$. By
rearrangement we obtain $\phi_{I'}\circ\phii_I=\l^{k'-k}$ and  $\l^{k'-k}\circ\phi_I=\phi_{I'}$.
%
\qed


Feng and Wang \cite[Theorem 1.1 (The Logarithmic Commensurability
Theorem)]{FW} proved that if $K=\phi_1(K)\cup\ldots\cup\phi_r(K)$ is a
self-similar set in $\R$ satisfying the open set condition with 
Hausdorff dimension less than $1$ and such that each similarity
map $\phi_i$ is of the form $\phi_i(x)=bx+c_i$ with a fixed $b$ 
and $a K + t\su K$ for some $a,t\in\R$ then $\log|a|/\log|b|\in\Q$.
They also posed the problem (Open Question 2) of generalizing this result to
higher dimensions. If we assume the strong open set condition instead of the 
open set condition then the above Theorem~\ref{lambdatul} tells much more
about the maps $\phi_1,\ldots,\phi_r$ and $ax+t$ and immediately gives
the following higher dimensional generalization of the
Logarithmic Commensurability
Theorem of Feng and Wang, in which 
we can also allow non-homogeneous self-similar sets. 

\begin{cccc}\label{log}
Let $\egy$ be a self-similar set with strong separation condition and 
suppose that $\l$ is a similitude for which $\l(K)\subset K$.
If $a_1,\ldots, a_r$ and $b$ denote the similarity ratios of 
$\phi_1, \ldots, \phi_r$ and $\la$, respectively, then $\log b$ must be
a linear combination of $\log a_1,\ldots, \log a_r$ with rational coefficients.
\end{cccc}


\section{Isometry invariant measures}\label{measures}

In this section all self-similar sets we consider will satisfy the strong separation condition, even if we do not mention it every time.

Before we start to study and characterize the isometry invariant measures on a self-similar set of strong separation condition, we have to pay some attention to the connection of a self-similar set and the self-similar measures living on it.

We have called a compact set $K$ self-similar with SSC if $\egy$ holds for some
similitudes $\phi_1, \ldots, \phi_r$. A \emph{presentation} of $K$ is a finite collection of
similitudes
$\{\psi_1, \ldots, \psi_s\}$, such that $\egys$ and $s\ge 2$. Clearly, a self-similar set with SSC has many different presentations. For example, if $\{\phik\}$ is a presentation of $K$, then $\{\phi_i\circ\phi_j: 1\le i,j \le r\}$ is also a presentation.

As we shall see in the next section, it can even happen that a self-similar set has no ``smallest'' presentation.
We say that a presentation $\F_1=\{\psik\}$ is smaller than the presentation $\F=\{\phik\}$, if for every $1\le i\le r$ there exists a multi-index $I$, such that $\phi_i=\psi_I$. This defines a partial ordering on the presentations; let us denote by $\F_1 \le \F$ if $\F_1$ is smaller than $\F$. We call a presentation \emph{minimal}, if there is no smaller presentation (excluding itself). We call a presentation \emph{smallest}, if it is smaller than any other presentation.

There exists a self-similar set with SSC which has more than one minimal presentations; that is, it has no smallest presentation (see Section~\ref{pelda}).

The notion of a self-similar measure on a self-similar set depends on the presentation. Thus, when we say that $\mu$ is a self-similar measure on $K$, we always mean that $\mu$ is self-similar measure with respect to the given presentation of $K$. Clearly if $\F_1\le \F$, then there are less self-similar measures with respect to $\F_1$ than to $\F$.
 It will turn out that the isometry invariant self-similar measures are the same independently of the presentations.


\begin{notation}\upshape
For the sake of simplicity, for a similitude $\l$ with $\l(K)\subset K$ let $\mu(\l)$ denote $\m{\l(K)}$.
In the composition of similitudes we might omit the mark $\circ$, so $g_1 g_2$ stands for $g_1\circ g_2$, and by $g^k$ we shall mean the composition of $k$ many $g$'s.
\end{notation}

Clearly, given any self-similar measure $\mu$, $\mu\circ\phi_I=\mu(\phi_I)\cdot\mu$ holds for the similitudes $\phi_I$ arising from the presentation of $K$.
According to the next proposition, if for a given self-similar measure $\mu$ the congruent elementary pieces are of equal measure, then the same holds for any similitude $\l$ satisfying $\l(K)\subset K$; that is, we have
$\mu\circ\l=\mu(\l)\cdot\mu$ as well.
%

\begin{all}\label{mong}
Let $\egy$ be self-similar set with strong separation condition, and $\mu$ be a self-similar measure on $K$ for which the congruent elementary pieces are of equal measure.
\begin{enumerate}
\item Then for every similitude $\l$ with $\l(K) \subset K$, $\,\mu \circ \l = \m{\l(K)}\cdot \mu$ holds; that is,
for any Borel set $H\subset K$ we have $\m{\l(H)} = \m{\l(K)}\cdot \m{H}$.
\item For every elementary piece $\phi_I(K)$ and for every isometry $g$ for which $g(\phi_I(K)) \subset K$ holds,
we have $\m{\phiIK}=\m{g(\phiIK)}$.
\end{enumerate}
\end{all}

\Biz
%
%
According to Lemma~\ref{vegessok} there are only finitely many similitudes $\l$ for which
$\l(K)\subset K$ holds and $\l(K)$ intersects at least two first generation elementary pieces.
Denote these by $\l_0, \l_1, \ldots, \l_t$, where $\l_0$ should stand for the identity. 

We claim that it is enough to prove the first part of the proposition for these similitudes only. Let $\l$ be a
similitude for which $\l(K)\subset K$. Let $\phi_I(K)$ the smallest elementary piece which contains $\l(K)$.
Then the similitude $\phi_I^{-1}\circ\l$ maps $K$ into itself and the image intersects at least two first generation elementary pieces, hence it is equal to a similitude $\l_i$ for some $i$. Thus $\l=\phi_I\circ\l_i$.
The measure $\mu$ being self-similar we have $\mu\circ\phi_J=p_J\cdot\mu=\m{\phi_J(K)}\cdot\mu$ for every multi-index $J$, hence for any Borel set $H\subset K$ we obtain
\mult{
\m{\l(H)} = \m{(\phi_I \circ \l_i)(H)}=\m{\phi_I(K)} \cdot \m{\l_i(H)} = \m{\phi_I(K)} \cdot \m{\l_i(K)} \cdot \m{H} \\
= \m{(\phi_I\circ\l_i)(K)} \cdot \m{H} = \m{\l(K)}\cdot\m{H},
}
as we stated.

According to Theorem~\ref{lambdatul} for every integer $i$ with $0\le i \le t$ there exist multi-indices $I_i$, $J_i$ and a positive integer $k_i$, for which $\l_i^{k_i} \circ \phiIi=\phiJi$. Let $b_i\=\phi_{I_i}$, $c_i\=\phi_{J_i}$, hence $\l_i^{k_i} b_i =c_i$.

Let $$\mu^*(\l_i) \= \sqrt[k_i]{\frac{\mu(c_i)}{\mu(b_i)}}.$$
Our aim is to show that $\mu^*(\l_i)=\mu(\l_i)$.

For every integer $i$ with $0\le i \le t$ and for every multi-index $I$ there exists an integer $j$, $0 \le j \le t$, and a multi-index $J$ such that $\l_i \circ \phi_I = \phi_J \circ \l_j$ (let $\phi_J(K)$ be the smallest elementary piece which contains $(\l_i \circ \phi_I)(K)$).

We define the \emph{congruency} equivalence relation among similitudes: for similitudes $g_1$ and $g_2$
let $g_1 \kb g_2$ denote that $g_1 \circ g_2^{-1}$ is an isometry; that is, for every set $H$ the sets $g_1(H)$ and
$g_2(H)$ are congruent. This is the same as that the similarity ratio of $g_1$ and $g_2$ are equal. Hence congruency
is independent of the order of the composition, so $g_1\circ g_2 \kb g_3 \Longleftrightarrow g_2\circ g_1 \kb g_3$.
Using the equalities $\l_i \phi_I = \phi_J \l_j$, $\l_i^{k_i}b_i=c_i$ and $\l_j^{k_j}b_j=c_j$ we obtain
$$\underbrace{ \l_i^{k_i k_j} \,\phi_I^{k_i k_j} }_{\kb \phi_J^{k_i k_j} \l_j^{k_i k_j}} \, b_i^{k_j} \, b_j^{k_i} \ \kb \ \phi_J^{k_i k_j} \, b_i^{k_j} \, \underbrace{\l_j^{k_i k_j} \, b_j^{k_i}}_{\kb c_j^{k_i}} \ \kb \ \phi_J^{k_i k_j} \, b_i^{k_j} \, c_j^{k_i},$$
%
$$\l_i^{k_i k_j} \,\phi_I^{k_i k_j} \,b_i^{k_j} \, b_j^{k_i} \ \kb\ 
\underbrace{\l_i^{k_i k_j} \, b_i^{k_j}}_{\kb c_i^{k_j}} \, \phi_I^{k_i k_j} \, b_j^{k_i} \ \kb \ c_i^{k_j} \, \phi_I^{k_i k_j} \, b_j^{k_i}.$$
Comparing these we get
$$\phi_J^{k_i k_j} \, b_i^{k_j} \, c_j^{k_i} \ \kb \ c_i^{k_j} \, \phi_I^{k_i k_j} \, b_j^{k_i}.$$
Since all the similitudes $b_i$, $b_j$, $c_i$, $c_j$ are some composition of similitudes of the presentation,
%
the elementary pieces $\big(\phi_J^{k_i k_j} \, b_i^{k_j} \, c_j^{k_i}\big)(K)$ and $\big(c_i^{k_j} \, \phi_I^{k_i k_j} \, b_j^{k_i}\big)(K)$ are congruent, so they are of equal measure. The measure is self-similar, thus
%
%
$$\mm{\phi_J}^{k_i k_j} \mm{b_i}^{k_j} \mm{c_j}^{k_i} = \mm{c_i}^{k_j} \mm{\phi_I}^{k_i k_j} \mm{b_j}^{k_i},$$
hence by the definition of $\mu^{*}$ we get
%
$$\mm{\phi_J}^{k_i k_j} \mmcs{\l_j}^{k_i k_j} = \mmcs{\l_i}^{k_i k_j} \mm{\phi_I}^{k_i k_j},$$
$$\mmcs{\l_j} \mm{\phi_J} = \mmcs{\l_i} \mm{\phi_I}.$$
Therefore
$$\mm{\l_i \phi_I} = \mm{\phi_J \l_j} = \mm{\phi_J} \mm{\l_j} = \frac{\mmcs{\l_i} \m{\phi_I}}{\mcs{\l_j}} \m{\l_j}.$$
Altering this we got the following: for every $i$ and $I$ there exists $j$ such that
%
$$\mm{\l_i \phi_I} = \mmcs{\l_i} \frac{\mm{\l_j}}{\mmcs{\l_j}} \mm{\phi_I}.$$
Note that $\mmcs{\l_j}\neq 0$.

Let $m$ be an index for which $$\frac{\mm{\l_m}}{\mmcs{\l_m}} \lep \frac{\mm{\l_i}}{\mmcs{\l_i}}$$ for every index $0 \le i \le t$. We label some inequalities by a dot so we can refer to them later.
Then for any $\phi_I$,
$$\mm{\l_m \phi_I} = \mmcs{\l_m} \frac{\mm{\l_j}}{\mmcs{\l_j}} \mm{\phi_I} \gep \mmcs{\l_m} \frac{\mm{\l_m}}{\mmcs{\l_m}} \mm{\phi_I} = \mm{\l_m} \mm{\phi_I}$$
for some index $j$ with $0 \le j \le t$.

Let $\left\{\phiIiK \right\}$ be a finite partition of $K$ with elementary pieces such that the partition includes $\phi_I(K)$. Then
\mult{
\m{\l_m(K)}=\M{\l_m\left(\bigcupp \phiIiK\right)}= \M{\bigcupp \l_m(\phiIiK)}=
\sum \mm{\l_m \phiIi} \\
\gep \sum \mm{\l_m} \mm{\phiIi}=\mm{\l_m},
}
hence equality holds everywhere, so $\mm{\l_m \phi_I} = \mm{\l_m} \mm{\phi_I}$ for every multi-index $I$.

Let $H\subset K$ be a Borel set. By the definition of the measure $\mu$, there exist elementary pieces $a_{ij}(K)$
for which $H\subset \bigcap_j \bigcup^{*}_i a_{ij}(K)$ and $\mm{H} = \inf_j \m{\bigcup^{*}_i a_{ij}(K)}= \m{\bigcap_j \bigcup^{*}_i a_{ij}(K)}$ hold. Then
\mult{
\mm{\lm(H)} \le 
\M{\lm\Big(\bigcap_j \bigcupp_i a_{ij}(K)\Big)} = 
\M{\bigcap_j \bigcupp_i \lm(a_{ij}(K))} \\
\le \inf_j \M{\bigcupp_i \lm(a_{ij}(K))}
= \inf_j \sum_i \mm{\lm a_{ij}} = \inf_j \sum_i \mm{\lm} \mm{a_{ij}} \\
= \mm{\lm} \inf_j \sum_i \mm{a_{ij}}= \mm{\lm} \M{\bigcap_j \bigcupp_i a_{ij}(K)}=
\mm{\lm} \mm{H}.
}
Repeating this argument for $H^c\=K\bol H$ we obtain $\m{\lm(H^c)} \le \mm{\lm} \mm{H^c}$.
Summing these we get $\m{\lm(H))} + \m{\lm(H^c)} \le \mm{\lm} \mm{H} + \mm{\lm} \mm{H^c}$,
in fact this is an equality, so we have $\m{\lm(H)}=\mm{\lm} \mm{H}$. Thus $\mu\circ\lm=\mm{\lm}\cdot \mu$.

From this we obtain that for any Borel set $H\subset K$,
$$\m{\lm^n(H)}=\m{\lm(\lm^{n-1}(H))}=\mm{\lm} \m{\lm^{n-1}(H)},$$
and by induction we get that $\m{\lm^n(H)}=\mm{\lm}^n \mm{H}$, hence $\mm{\lm^n}=\mm{\lm}^n$.

Therefore $\mm{\lm^{k_m} b_m} = \mm{\lm}^{k_m} \mm{b_m}$ holds. From the definition of $\mmcs{\lm}$ we have 
$\mm{c_m}=\mmcs{\lm}^{k_m} \mm{b_m}$ and $c_m=\lm^{k_m} b_m$, thus
$$\mm{\lm}^{k_m} \mm{b_m} = \mm{\lm^{k_m} b_m} = \mm{c_m} = \mmcs{\lm}^{k_m} \mm{b_m}.$$
Since $\mm{b_m}>0$, we get $\mm{\lm}=\mmcs{\lm}$. Since $m$ was chosen to be that index $i$ for which $\frac{\mm{\l_i}}{\mmcs{\l_i}}$ is minimal, we get that $\mmcs{\l_i} \lep \mm{\l_i}$ for every $0\le i \le t$.
%

Now we can repeat the whole argument for such an index $m$ for which $\frac{\mm{\l_m}}{\mmcs{\l_m}} \ge \frac{\mm{\l_i}}{\mmcs{\l_i}}$ holds for every index $i$ ($0 \le i \le t$). We just have to reverse the inequalities labelled with a dot, and we obtain that for every index $i$ ($0\le i \le t$), $\mmcs{\l_i} \ge \mm{\l_i}$ holds.
Thus for every $i$ ($0\le i \le t$) we have $\mmcs{\l_i} = \mm{\l_i}$. Therefore we could choose any $i$ ($0\le i \le t$) as $m$, so for every $i$ the equality $\mu \circ\l_i = \mu(\l_i) \cdot \mu$ holds. By the observation we made at the beginning of the proof we get that for every
similitude $\l$ with $\l(K)\subset K$, $\mu\circ\l=\mu(\l)\cdot\mu$ holds, thus $\mu\circ\l^n=\mu(\l)^n\cdot\mu$ holds as well for any positive integer $n$.


Now we shall prove the second part of the proposition. Suppose that the isometry $g$ maps the elementary piece $\phi_L(K)$ into $K$, so $g(\phi_L(K))\subset K$. By Theorem~\ref{lambdatul} there exist multi-indices $I$, $J$ and a positive integer $k$ such that $(g\circ\phi_L)^k \circ \phi_I=\phi_J$. Using the first part of this proposition (which is already proven) we get
%
\begin{equation}\label{a}
\mu(\phi_J)=\m{(g\circ\phi_L)^k \circ \phi_I}=\mm{g\circ\phi_L}^k \mm{\phi_I}.
\end{equation}
Clearly $\phi_J=(g\circ\phi_L)^k \circ \phi_I \kb (\phi_L)^k \phi_I$, thus
\begin{equation}\label{b}
\mu(\phi_J)= \m{(\phi_L)^k \phi_I}=\mu(\phi_L)^k \mu(\phi_I).
\end{equation}
By \eqref{a} and \eqref{b} we obtain
$$\mm{g\circ\phi_L}^k \mm{\phi_I} = \mu(\phi_L)^k \mu(\phi_I),$$
$$\mm{g\circ\phi_L}=\mu(\phi_L),$$
which proves the proposition.
\qed

\begin{tetel}[Characterization of isometry invariant measures]\label{karakterizacio}
\ \\
Let
$\egy$ be a self-similar set with the strong separation condition and $\mu$ a self-similar measure on $K$ for which congruent elementary pieces are of equal measure. Then $\mu$ is an isometry invariant measure on $K$.
\end{tetel}
\Biz
We have to show that for any isometry $g$ and Borel set $H\subset K$ if $g(H)\subset K$ then $\mu(H)=\m{g(H)}$.

Let $c<1$ be the constant given by Theorem~\ref{andras}. At first consider a set $H\subset K$ of positive measure.
Applying Lemma~\ref{suruseg} for the set $H$ with $\eps=1-c$ we obtain that there exists an elementary piece $a(K)$ for which $\m{H\cap a(K)}>c\cdot\m{a(K)}$.
Since $H\subset g^{-1}(K)$, we have $\m{g^{-1}(K) \cap a(K)}>c\cdot\m{a(K)}$, so applying Theorem~\ref{andras} $a(K) \subset g^{-1}(K)$, $g(a(K))\subset K$. Put $\l = g\circ a$. According to the second part of Proposition~\ref{mong} we have $\mu(\l)=\mu(a)$ (where $\mu(\l)$ is an abbreviation of $\m{\l(K)}$), and putting $H_0\=a^{-1}(a(K)\cap H)$ we have $\m{\l(H_0)}=\mm{\l}\mm{H_0}$, thus
\mult{
0<c\cdot\m{a(K)}<\m{a(K)\cap H}=\m{a(H_0)}=\mm{a}\mm{H_0}=\mm{\l}\mm{H_0}\\
= \m{\l(H_0)}= \m{g(a(H_0))}=\m{g(a(K) \cap H)} \le \m{g(H)},
}
so $g(H)$ is of positive measure. Thus a congruent copy of a set of positive measure is of positive measure, and a congruent copy of a negligible set is also negligible.

Now let $H\subset K$ be any Borel set, $g$ an isometry, for which $g(H)\subset K$. Apply Lemma~\ref{moho} with some $0<\eps<1-c$. We obtain elementary pieces $a_i(K)$ such that
%
$$\m{H \cap a_i(K)} > (1-\eps)\cdot \m{a_i(K)} \quad \text{and} \quad \m{H\bol \bigcupp_i a_i(K)}=0.$$
%
Then $H\subset g^{-1}(K)$, therefore $\m{g^{-1}(K) \cap a_i(K)}>(1-\eps)\cdot\m{a_i(K)}$. 
According to Theorem~\ref{andras},
$g^{-1}(K) \supset a_i(K)$, so $g(a_i(K))\subset K$. By the second part of Proposition~\ref{mong} we get $\m{g(a_i(K))}=\m{a_i(K)}$, and using the fact that a congruent copy of a set of zero measure is also of zero measure,
\mult{
\m{g(H)}=\M{g\Big(H\cap\bigcupp_i\! a_i(K)\Big)} + \M{g\Big(H\bol \bigcupp_i\!\! a_i(K)\Big)}=
\M{g\Big(H\cap\bigcupp_i \! a_i(K)\Big)} \\
= \sum_i \m{g(H\cap a_i(K))} \le \sum_i \m{g(a_i(K))}
= \sum_i \m{a_i(K)} \\
\le \frac{1}{1-\eps}\cdot \sum_i \m{H \cap a_i(K)}= \frac{1}{1-\eps}\cdot
\M{H\cap \bigcupp_i a_i(K)} = \frac{1}{1-\eps}\cdot \mm{H}.
}
This is true for any $0<\eps<1-c$, hence $\m{g(H)} \le \m{H}$. Repeating this argument for $g(H)$ instead of $H$ and for $g^{-1}$ instead of $g$ gives $\m{H} \le \m{g(H)}$, hence $\m{H}=\m{g(H)}$. Thus $\mu$ is isometry invariant.
\qed

\br\label{magyarazas}
Using this theorem it is relatively easy to decide whether a self-similar measure is isometry invariant or not.
Denote the similarity ratio of the similitude $\phi_i$ by $\alfa_i$. It is clear that two elementary pieces are congruent if and only if they are image of $K$ by similitudes of equal similarity ratio.
Thus a self-similar measure $\mu$ is isometry invariant if and only if
provided that
$\alfa_{i_1} \alfa_{i_2} \ldots \alfa_{i_n}= \alfa_{j_1} \alfa_{j_2} \ldots \alfa_{j_m}$ holds,
the equality  $p_{i_1} p_{i_2} \ldots p_{i_n}= p_{j_1} p_{j_2} \ldots p_{j_m}$ also holds for the weights of the measure $\mu$.
By switching from the similarity ratios $\alfa_i$ and weights $p_i$ to the negative of their logarithm we get a system of linear equations for the variables $-\log p_i$.
The solutions of this system (which also satisfy the normalizing equation $\sum_i p_i=1$) give those weight vectors which define isometry invariant measures on $K$. 

For example, it is easy to see that if the positive numbers $-\log \alfa_i$ $(i=1,\ldots, r)$ are linearly
independent over $\mathbb{Q}$, then every self-similar measure is isometry invariant.


So, to the $r$ dimensional vectors, formed by the $-\log p_i$ weights of the isometry invariant measures, correspond the intersection of a linear subspace of $\R^r$ and the hypersurface corresponding to $\sum_i p_i= 1$. That this subspace is of dimension at least $1$ and intersects the positive part of the space $\R^r$, we know from the existence of Hausdorff measure. (Or rather from the fact that the weights $p_i=\alfa_i^s$ automatically satisfy all the equalities.)
\er

The notion of a self-similar measure depended on the the choice of the presentation. However, the converse is true for the notion of an isometry invariant self-similar measure.

\begin{tetel}\label{vanertelme}
Let $K$ be self-similar with the strong separation condition and $\{\phik\}$ a presentation of it. Let $\mu$ be isometry invariant and self-similar with respect to this presentation. Then $\mu$ is self-similar with respect to any presentation of $K$. Thus the class of isometry invariant self-similar measures is independent of the choice of presentation.
\end{tetel}
\Biz
Let $\{\psi_1, \ldots, \psi_s\}$ be an other presentation of $K$.
According to Theorem~\ref{lambdatul} there exist positive integer $k$ and elementary pieces $\phi_I$, $\phi_J$
such that $\psi_i^k \circ \phi_I = \phi_J$, so applying the first part of Proposition~\ref{mong} we get
%
$$0<\mu(\phi_J) = \mu(\psi_i^k \circ \phi_I) = \mu(\psi_i)^k \mu(\phi_I),$$
that is, $\mu(\psi_i)>0$ for every $1\le i\le s$.

According to the first part of Proposition~\ref{mong}, $\mu \circ \psi_i = \mu(\psi_i) \cdot \mu$, and since $\sum \m{\psi_i(K)}=1$ holds, this means exactly that $\mu$ is a self-similar measure with respect to the presentation $\{\psi_1, \ldots, \psi_s\}$.
\qed

\begin{ddd}\upshape
Let $\egy$ be a self-similar set with strong separation condition.
Put $S=\{-\log \alpha_i:1\le i\le r\}$, where $\alpha_i$ is the similarity ratio of $\phi_i$.
The \emph{algebraic dependence number} (of this presentation) is the dimension over $\Q$ of the vectorspace
generated by $S$ minus one.
\end{ddd}

By Remark~\ref{magyarazas} it is easy to see that the algebraic dependence number of a presentation is exactly the same as the topological dimension of the surface corresponding to the isometry invariant self-similar measures on $K$.
Thus, by Theorem~\ref{vanertelme}, one can prove the following.

\begin{theorem}\label{algebraic}
The algebraic dependence number of a self-similar set does not depend on the presentation we choose.
\end{theorem}

We mention that it is easy to show that the algebraic dependence number is the same for two presentations $\F_1$ and $\F_2$ if $\F_1\le \F_2$; that is, when one of them extends the other in the trivial way we defined at the beginning of this section. However, there are self-similar sets with two presentations which have no common extension and they are not an extension of the same third presentation (see Theorem~\ref{remelemigaz}). Thus we have no direct (or trivial) proof for Theorem~\ref{algebraic}.

An easy consequence of the characterization theorem is the following.

\begin{kov}\label{iranyitastarto}
Let $\egy$ be a self-similar set with strong separation condition, $\mu$ be a self-similar measure on $K$. Then if $\mu$ is invariant under orientation preserving isometries, then it is invariant under all isometries.
\end{kov}
\Biz
According to Theorem~\ref{karakterizacio} it is enough to show that congruent elementary pieces are of equal measure. Let $\phi_I(K)$ and $\phi_J(K)$ be congruent elementary pieces. Then $\phi_I^2(K)$ and $\phi_J^2(K)$ are also congruent elementary pieces, $\phi_I^2$ and $\phi_J^2$ are orientation preserving similitudes, so $\phi_I^2 \circ \phi_J^{-2}$ is an orientation preserving isometry, hence by the assumption $\m{\phi_I^2(K)}=\m{\phi_J^2(K)}$.
Since $\mu$ is self-similar, $\m{\phi_I^2(K)}=\m{\phi_I(K)}^2$ and $\m{\phi_J^2(K)}=\m{\phi_J(K)}^2$, thus $\m{\phi_I(K)} = \m{\phi_J(K)}$. This proves the statement.
\qed

\section{Minimal presentations}\label{pelda}

At first we give an example for a self-similar set on the line (with strong separation condition) which has no smallest presentation, that is, it has more than one minimal presentations.
Set $\phi_1(x)=\frac{x}{3}$, $\phi_2(x)=\frac{x}{3}+\frac{2}{3}$, let $K$ be the compact set for which $K=\phi_1(K)\cup \phi_2(K)$, apparently this is the triadic Cantor set. Set $\psi_1(x)=-\frac{x}{3}+\frac{1}{3}$. Then $K=\psi_1(K)\cupp \phi_2(K)$ as well, and it is clear, that both of these two different presentations are minimal, since they consist of only two similitudes.

However, these two presentations are not ``essentially different'': the sets $\{\phi_1(K), \phi_2(K)\}$ and $\{\psi_1(K), \phi_2(K)\}$ coincide. On \emph{essential presentation} we shall mean not the set of the similitudes but rather the set of the first generation elementary pieces. We shall say that the essential presentation $\{a_1(K), \ldots, a_r(K)\}$ is \emph{briefer} than the essential presentation $\{b_1(K), \ldots, b_s(K)\}$, if for every $j=1, \ldots, s$ there exists $1\le i\le r$ such that $b_j(K)\subset a_i(K)$. We call an essential presentation minimal if the only briefer essential presentation is itself, and we call it the smallest if it is briefer than any other essential presentation. It is easy to check that the triadic Cantor set possesses a smallest essential presentation.

In the followings we shall present a self-similar set which has got no smallest essential presentation, that is, it has minimal essential presentations more than one.

\br
The following statement is true for many self-similar sets $K$: If $\l_1$ and $\l_2$ are similitudes for which
$\l_1(K)\subset K$, $\l_2(K)\subset K$ and $\l_1(K)\cap \l_2(K) \neq \emptyset$,
then $\l_1(K)\subset\l_2(K)$ or $\l_2(K)\subset\l_1(K)$.
The proofs of Section~\ref{similar} would have been much simpler if this statement has held for every self-similar set satisfying the strong separation condition. However this statement does not hold generally as we shall show in our following construction. We note that this statement is not necessarily equivalent to that $K$ has only one minimal essential presentation. See also the end of Section~\ref{conc} and especially Question~\ref{relnyilt2}.
\er

\begin{theorem}
There exists a self-similar set $K$ with the strong separation condition which has no smallest essential presentation. Moreover, there exists similitudes $\l_1$ and $\l_2$ such that
$\l_1(K)\cap \l_2(K) \neq \emptyset$, but
$\l_1(K)\not\subset\l_2(K)$ and $\l_2(K)\not\subset\l_1(K)$.
\end{theorem}

\bp
We present a figure of our construction. One may check the proof of this theorem just by looking at that figure.

Let $a,b,c$ positive integers for which $a+b+a+c+a+b+a=1$ \,and \,$b=a\cdot c$. It is easy to see that for every $0<a<1/4$ there exist a unique $b$ and $c$ with these conditions.
Let $\phi_1$ be the orientation preserving similitude mapping the interval $[0,1]$ onto the interval $[0,a]$.
Let $\phi_2$ take the interval $[0,1]$ onto $[a+b,a+b+a]$, $\phi_3$ onto $[1-a-b-a, 1-a-b]$, and $\phi_4$ onto $[1-a,1]$, all of them preserving the orientation.
That is, $\phi_1(x)=a\cdot x$, $\phi_2(x)=a\cdot x+ a+b$, $\phi_3(x)=a\cdot x+1-a-b-a$, $\phi_4(x)=a\cdot x+1-a$.



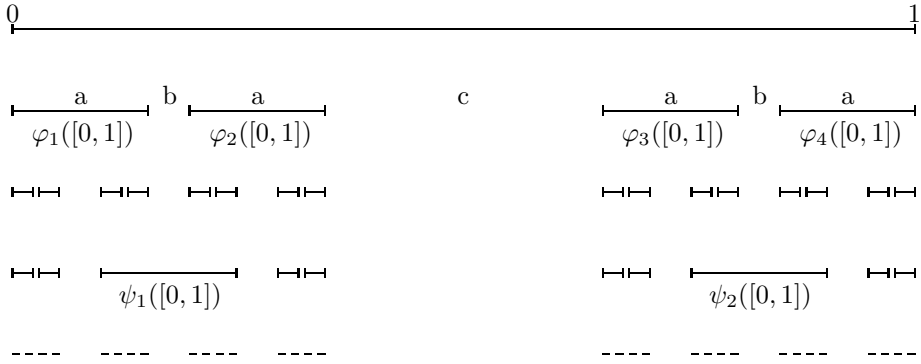
\setlength{\unitlength}{12cm}
\begin{figure}[h] \centering
\begin{picture}(0,+.4200)(0,-0.3000)

\input{rajz}
\put(-.5070,+0.098){0}
\put(+.4920,+0.098){1}
\def\fel{0.007}
\put(-.43195,\fel){a}
\put(-.33387,\fel){b}
\put(-.2358,\fel){a}
\put(-.00695,\fel){c}
\put(.41805,\fel){a}
\put(.31996,\fel){b}
\put(.22189,\fel){a}

\put(-.4788,-0.034){$\phi_1([0,1])$}
\put(-.282,-0.034){$\phi_2([0,1])$}
\put(.1750,-0.034){$\phi_3([0,1])$}
\put(.3718,-0.034){$\phi_4([0,1])$}

\put(-.382,-.214){$\psi_1([0,1])$}
\put(.2718,-.214){$\psi_2([0,1])$}

\end{picture}
\end{figure}

Let $K$ be the unique compact set for which $K=\phi_1(K) \cupp \phi_2(K) \cupp \phi_3(K) \cupp \phi_4(K)$. Thus the first generation elementary pieces of $K$ are of diameter $a$, and there are ,,holes'' between them of length $b$, $c$ and $b$. 
It is clear that $K\subset [0,1]$ and $K$ is symmetric to $\frac{1}{2}$.

The second row of the figure symbolizes this presentation of $K$, more precisely it shows the intervals $\phi_i([0,1])$ (choosing $a=0.15$, $\,c=\frac{0.4}{1.3}$).
In the first row the interval $[0,1]$ can be seen. The third row of the figure shows the intervals $\phi_i(\phi_j([0,1])) \quad (1\le i,j \le 4)$. The fifth row tries to present the set $K$.

Set $\psi_1(x)=a\cdot x +a^2+a\cdot b+a^2+a\cdot c$ and $\psi_2(x)=a\cdot x +1-a-b-a^2-a\cdot b -a^2$. In the fourth row of the figure the images of the interval $[0,1]$ by the similitudes
$\phi_1^2$, $\phi_1\circ\phi_2$, $\psi_1$, $\phi_2\circ\phi_3$, $\phi_2 \circ \phi_4$, $\phi_3\circ\phi_1$, $\phi_3 \circ \phi_2$, $\psi_2$, $\phi_4 \circ \phi_3$, $\phi_4^2$
are shown.

We claim that $\psi_1(K)\subset K$ and $\psi_2(K) \subset K$, moreover $$\{\phi_1^2,  \,\phi_1\circ\phi_2, \,\psi_1, \,\phi_2\circ\phi_3, \,\phi_2 \circ \phi_4, \,\phi_3\circ\phi_1, \,\phi_3 \circ \phi_2, \,\psi_2, \,\phi_4 \circ \phi_3, \,\phi_4^2\}$$
is a presentation of $K$ (see the fourth row of the figure). For this it is sufficient to prove that 
$\psi_1 \circ \phi_1= \phi_1 \circ \phi_3$,
$\psi_1 \circ \phi_2= \phi_1 \circ \phi_4$,
$\psi_1 \circ \phi_3= \phi_2 \circ \phi_1$,
$\psi_1 \circ \phi_4= \phi_2 \circ \phi_2$,
and
$\psi_2 \circ \phi_1= \phi_3 \circ \phi_3$,
$\psi_2 \circ \phi_2= \phi_3 \circ \phi_4$,
$\psi_2 \circ \phi_3= \phi_4 \circ \phi_1$,
$\psi_2 \circ \phi_4= \phi_4 \circ \phi_2$.
These can be easily checked, all equalities rely on the choice of $b=a\cdot c$.

Now we prove that there does not exist an essential presentation $\{\varrho_1(K), \ldots, \varrho_r(K)\}$ of the self-similar set $K$ which is briefer than both of the essential presentations corresponding to the original
presentation $\{\phi_1, \phi_2, \phi_3, \phi_4\}$ and the presentation just defined above. This would prove that $K$ has no unique minimal essential presentation. (In fact both of these essential presentations are minimal.) 
Indirectly suppose that there exists an essential presentation $\{\varrho_1(K), \ldots, \varrho_r(K)\}$ of this kind.
Since $\phi_1(K)\cap \psi_1(K)\neq \ures$, for some $i$ $\phi_1(K)\cup \psi_1(K) \subset \varrho_i(K)$.
For the same $i$ we also have $\phi_2(K)\cup \psi_1(K) \subset \varrho_i(K)$.
Similarly there exists an index $j$ such that $\phi_3(K) \cup \psi_2(K) \cup \phi_4(K) \subset \varrho_j(K)$. From this we conclude that $K=\varrho_i(K) \cupp \varrho_j(K)$, but then the similitudes $\varrho_i$ and $\varrho_j$ could only be the ones mapping $[0,1]$ onto $[0, a+b+a]$ and $[1-a-b-a,1]$. This yields to $b=(a+b+a)\cdot c$, which contradicts $b=a\cdot c$.

The similitudes $\l_1$ and $\l_2$ we promised can chosen to be $\phi_1$ and $\psi_1$.
\qed
\br
This example (and many other results of the present article) is contained in
the Master Thesis of the third author \cite{M}. Independently, Feng and Wang
in \cite{FW} exhibit an almost identical example. Moreover, much of their
paper is devoted to the investigation of the structure of possible
presentations of given self-similar sets; or, using their terminology, 
the structure
of generating iterated function systems of self-similar sets. They also prove
positive results (that is, when a smallest presentation does exist) under
various assumptions. 
\er

\begin{theorem}\label{remelemigaz}
There exists a self-similar set $K$ with the strong separation condition and two (essential) presentations of $K$, $\F_1$ and $\F_2$,
such that there is no presentation $\G$ which is a common extension of $\F_1$ and $\F_2$,
nor there exists an (essential) presentation which is smaller (briefer) than $\F_1$ and $\F_2$.
\end{theorem}

Thus, Theorem~\ref{algebraic} cannot be proved in the trivial way (see our remarks after that theorem).
We leave the proof of Theorem~\ref{remelemigaz} to the reader, with the instructions that one should choose the self-similar set $K$ constructed above, and the presentations of the second and fourth row of the figure should be chosen as $\F_1$ and $\F_2$.


\section{Intersection of translates of a self-affine Sierpi\'nski sponge}
\label{intersection}

The following is the key lemma for all results of this section.

\begin{propo}
\label{prop}
Let $K=K(M,D)$ and $\mu$ be like in Definition~\ref{def:K} 
and let $t\in\R^n$ be 
such that  $\norma{M^k t} >0$ for every $k=1,2,\ldots$.

Then $\m{K\cap (K+t)}>0$ implies that there exists a 
$$
w\in\{-1,0,1\} \times \ldots \times \{-1,0,1\} \sm \{(0,\ldots,0)\}
$$
such that $D+w=D$ modulo $(m_1,\ldots,m_n)$; that is,
$$
D+w+M(\Z^n)=D+M(\Z^n)=D-w+M(\Z^n).
$$
\end{propo}

\bp
To make the argument intuitive and precise we shall present
the same proof in an informal and in a formal way separately.

\emph{The informal proof:}
According to Lemma~\ref{approx} and Lemma~\ref{suruseg} we
can find a $k$ such that 
$M^k t$ is not very close to any point of $\Z^n$,
and a $k-1$-th generation elementary part $S$
of $K$ in which the density of $K+t$ is almost $1$.
Then in all the $r$ $k$-th generation elementary parts of $K$
that are in $S$ the density of $K+t$ is still very close to $1$.

Each of these subparts intersect some $k$-th generation
elementary parts of $K+t$. The key observation is that there
are at most $2^n$ possible ways how these parts can intersect
each other.

Since $M^k t$ is not very close to the lattice points, these
intersections are intersections of sets similar to $K$ such
that one is always a not very close translate of the other.
Hence Lemma~\ref{continuity} implies that they cannot have big
intersection.

Since the density of $K+t$ is very close to $1$ in all 
$k$-th generation elementary parts of $K$ that are in $S$, this
implies that in the two directions for which the possible 
intersection has biggest measure, 
$K+t$ must have a $k$-th generation elementary part.

Hence we get two periods of the pattern $D$ such that their
difference $w$ is in $\{-1,0,1\}\times\ldots\times\{-1,0,1\}$.

\emph{The formal proof:}
Applying Lemma~\ref{continuity} for $\eps=1/(2\max(m_1,\ldots,m_n))$
we get a $0<\de<1$ such that 
\begin{equation}
\label{folyt}
\m{K\cap(K+u)}\le 1-\de \quad 
\textrm{whenever } |u|\ge\frac{1}{2 \max(m_1,\ldots,m_n)}.
\end{equation}
Applying Lemma~\ref{suruseg} for $B=(K+t)\cap K$ 
and $\eps=\frac{\de}{2^n r}$ and Lemma~\ref{approx} we get a $k\in\N$
and a $k-1$-th generation elementary part $S$ of $K$ such that
\begin{equation}
\label{density}
\m{S\cap(K+t)} > \frac{1-\frac{\de}{2^n r}}{r^{k-1}}
\end{equation}
and
\begin{equation}
\label{kozelito}
\norma{M^k t}>\frac{1}{2 \max(m_1,\ldots,m_n)}. 
\end{equation}

Let $\Phi$ be the similarity map which maps $S$ to $M(K)=K+D$; 
that is,
\begin{align*}
\Phi(x)&= M^k(x-(M^{-(k-1)}\al_{k-1}+\ldots+M^{-1}\al_1)) \\
&= M^k x - (M \al_{k-1}+ M^2 \al_{k-2} + \ldots + M^{k-1} \al_1),
\end{align*}
where $S=M^{k-1}(K)+M^{-(k-1)}\al_{k-1}+\ldots+M^{-1}\al_1$.

Using that $\Phi(S)=K+D=\cup_{j=1}^r K+d_j$, 
applying (\ref{mut}) and 
(\ref{density}) we get
\begin{align*}
  \mut\Big(\bigcup_{j=1}^r(K+d_j)\cap(\Phi(K+t))\Big)&=
   \mut\big(\Phi(S\cap(K+t))\big) = r^k \m{S\cap(K+t)} \\
   &> r^k\frac{1-\frac{\de}{2^n r}}{r^{k-1}} = r - \frac{\de}{2^n}.
\end{align*}

Since $\mut(K+d_j)=1$ $(j=1,\ldots,r)$ and the sets can intersect
each other only at a set of $\mut$-measure zero this implies that
\begin{equation}
\label{eq:8}
  \mut\big((K+d_j)\cap \Phi(K+t)\big) > 1 - \frac{\de}{2^n} \quad
\textrm{for every } j=1,\ldots,r.
\end{equation}

Since 
$\Phi(K)=M^k(K) - (M \al_{k-1} + \ldots + M^{k-1} \al_1)$ and
$M^k(K)\su K+D+ M(\Z^n)$, we have 
$\Phi(K)\su K+D+ M(\Z^n)$, and so 
$\Phi(K+t)\su K+D+ \Phi(t)+M(\Z^n)$.
Thus
\begin{multline*}
(K+d_j)\cap\Phi(K+t)\\ \su 
(K+d_j)\cap \big(K+D+ \Phi(t)+M(\Z^n)\big) \\
= \bigcup_{i=1}^r \big(K\cap (K+d_i+\Phi(t)-d_j+M(\Z^n))\big)+d_j.
\end{multline*}
Combining this with (\ref{eq:8}) and (\ref{mut}) (for $l=0$) we get
\begin{eqnarray}
  \label{eq:9}
\lefteqn{ 1-\frac{\de}{2^n} < \mut\big((K+d_j)\cap \Phi(K+t)\big)}\nonumber\\ 
& & \qquad \le \sum_{i=1}^r\mut\big(\big(K\cap(K+d_i+\Phi(t)-d_j+M(\Z^n))\big)+d_j\big)\\
& & \qquad \qquad = \sum_{i=1}^r\m{K\cap(K+d_i+\Phi(t)-d_j+M(\Z^n))}
\qquad (j=1,\ldots,r). \nonumber
\end{eqnarray}

Clearly, we have $\m{K\cap(K+d_i+\Phi(t)-d_j+M(\Z^n))}=0$ whenever
$$d_i+\Phi(t)-d_j \not\in (-1,1)\times\ldots\times(-1,1)+M(\Z^n).$$
Hence there are at most $2^n$ vectors $v\in\Z^n$
such that $v+\Phi(t)\in (-1,1)\times\ldots\times(-1,1)$;
let these vectors be $v_1,v_2,\ldots,v_p$, ($p\le 2^n$).

Thus, by omitting some zero terms on the right-hand side of (\ref{eq:9}) 
we can rewrite (\ref{eq:9}) as
\begin{equation}
\label{eq:10}
1-\frac{\de}{2^n} < 
\sum_{l\ :\ (\exists i)\ d_i-d_j\in v_l+M(\Z^n)}
\m{K\cap (K+v_l+\Phi(t))} \qquad (j=1,\ldots,r).
\end{equation}
 
Let 
$$
\be_l=\m{K\cap (K+v_l+\Phi(t))} \qquad (l=1,\ldots,p).
$$ 
By rearranging $v_1,\ldots,v_{p}$ if necessary, we may assume
that 
\begin{equation}
\label{eq:11}
\be_1\ge\be_2\ge\ldots\ge \be_{p}.
\end{equation}
Since $v_l\in\Z^n$ and $K\su[0,1]^n$, the sets $K+v_l+\Phi(t)$
($l=1,\ldots,p$) are pairwise disjoint and clearly 
$K=\cup_{l=1}^{p} K\cap(K+v_l+\Phi(t))$, we get
\begin{equation}
  \label{eq:12}
  1=\mu(K)=\sum_{l=1}^{p} \be_l.
\end{equation}

Since, using (\ref{kozelito}), 
$\norma{M^k t}>\frac{1}{2 \max(m_1,\ldots,m_n)}$, we have
$|v_1+\Phi(t)|>\frac1{2 \max(m_1,\ldots,m_n)}$.
Thus, by (\ref{folyt}), 
\begin{equation}
  \label{eq:13}
  \be_1=\m{K\cap(K+v_1+\Phi(t))}\le1-\de.
\end{equation}

Clearly (\ref{eq:11}), (\ref{eq:12}) and (\ref{eq:13}) 
implies that 
$\be_1\ge \be_2\ge \frac{\de}{p-1}>\frac{\de}{2^n}$ and so
\begin{align*}
  \be_1+\be_3+\be_4+\ldots \be_{2^n}&<
     1 - \frac{\de}{2^n} \qquad \textrm{ and}\\
  \be_2+\be_3+\be_4+\ldots \be_{2^n}&<
     1 - \frac{\de}{2^n}.
\end{align*}

Combining this with (\ref{eq:10}) we get that for every
$j\in\{1,\ldots,r\}$ there must be an $i_1$ such that
$d_{i_1}-d_j\in v_1+M(\Z^n)$ and an $i_2$ such that
$d_{i_2}-d_j\in v_2+M(\Z^n)$.
Since $D=\{d_1,\ldots,d_r\}$, this means that for every $d\in D$
we must have $d+v_1, d+v_2\in D+M(\Z^n)$.

Therefore $D+M(\Z^n) \supset D+v_1$ and so 
$D+M(\Z^n) \supset D+M(\Z^n) +v_1$. 
Applying this $m_1\cdot \ldots \cdot m_n$ many times we get
\begin{gather}
\label{eq:111}
D+M(\Z^n) \supset D+M(\Z^n) +v_1 \supset D+M(\Z^n) +2v_1 
\supset\ldots\nonumber \\
\ldots \supset D+M(\Z^n) 
+m_1\cdot \ldots \cdot m_n v_1 = D + M(\Z^n).
\end{gather}

Therefore $D+M(\Z^n)=D+M(\Z^n)+v_1$ and similarly
$D+M(\Z^n)=D+M(\Z^n)+v_2$. 
Thus $D+M(\Z^n)+v_1-v_2=D+M(\Z^n)=D+M(\Z^n)+v_2-v_1$.
Noting that, by definition, 
$w=v_1-v_2\in 
\{-1,0,1\}\times\ldots\times\{-1,0,1\}\sm\{(0,\ldots,0)\}$,
the proof is complete.
\qed

In order to use Proposition~\ref{prop} effectively we need
a discrete lemma.

\begin{lemma}
\label{discrete}
Let $M$ and $D$ be like in Definition~\ref{def:K}, 
$l\in\{1,2,\ldots,n\}$, $i\in\N$,
$$
D_i=M^{i-1}(D)+M^{i-2}(D)+\ldots+M(D)+D,
$$
and suppose that
\begin{equation}
\label{eq:112}
D_i+(\underbrace{ 1,\ldots,1 }_l,0,\ldots,0) +M^i(\Z^n)=D_i+M^i(\Z^n).
\end{equation}

Then at least one of the following two statements hold.

\noindent
(a) We have $m_1=\ldots=m_l$ and $a_1=\ldots=a_l$
for every $(a_1,\ldots,a_n)\in D$.

\noindent
(b) For some $l'\in\{1,2,\ldots,l-1\}$ we have 
\begin{equation*}
D_{i-1}+(\underbrace{ 1,\ldots,1 }_{l'},0,\ldots,0)+
M^{i-1}(\Z^n)=D_{i-1}+M^{i-1}(\Z^n).
\end{equation*}
\end{lemma}

\bp
Let $w=(\underbrace{ 1,\ldots,1 }_l,0,\ldots,0)$. From 
(\ref{eq:112}) we get
\begin{equation}
  \label{eq:113}
D_i+kw +M^i(\Z^n)=D_i+M^i(\Z^n) \qquad (k\in\Z).  
\end{equation}

First suppose that $a_1=\ldots=a_l$ does not hold for some
$a=(a_1,\ldots,a_n)\in D$. 
Then we can suppose that 
$a_1=\ldots=a_j<a_{j+1}\le\ldots\le a_l$ 
for some $j\in\{1,\ldots,l-1\}$.
Let $b=(b_1,\ldots,b_n)\in D_{i-1}$ be arbitrary.
Then $Mb+a\in M(D_{i-1})+D=D_i$.
Thus applying (\ref{eq:113}) for $k=-(a_1+1)$ we get
$$
Mb+a-(a_1+1)w \in D_i + M^i(\Z^n).
$$
Rewriting both sides we get
\begin{multline*}
M((b_1-1,\ldots,b_j-1,b_{j+1},\ldots,b_n)) \\
+(m_1-1,\ldots,m_j-1,a_{j+1}-a_1-1,
\ldots,a_l-a_1-1,a_{l+1},\ldots,a_n)\\
\in M(D_{i-1}+M^{i-1}(\Z^n))+D.
\end{multline*}
Since the second term of the left-hand side is in
$\{0,1\,\ldots,m_1-1\}\times\{0,1,\ldots,m_n-1\}$, 
we must have
$$
(b_1-1,\ldots,b_j-1,b_{j+1},\ldots,b_n)\in D_{i-1}+M^{i-1}(\Z^n).
$$
Since $b=(b_1,\ldots,b_n)\in D_{i-1}$ 
was arbitrary we get that
$$
D_{i-1}-(\underbrace{ 1,\ldots,1 }_{j},0,\ldots,0)
\su D_{i-1}+M^{i-1}(\Z^n),
$$
which implies, similarly like in (\ref{eq:111}), that 
$$
D_{i-1}+(\underbrace{ 1,\ldots,1 }_{j},0,\ldots,0)
+M^{i-1}(\Z^n) = D_{i-1}+M^{i-1}(\Z^n).
$$

Thus we proved that if $a_1=\ldots=a_l$ does not hold for
some $(a_1,\ldots,a_l)\in D$ then the statement (b) must hold. 
Exactly the same way (but ordering so that 
$m_1-a_1\le\ldots\le m_n-a_n$ and applying (\ref{eq:113}) 
for $k=m_1-a_1$ instead of $k=a_1$) we get that 
if $m_1-a_1=\ldots =m_l-a_l$ does not hold for some
$(a_1,\ldots,a_n)\in D$ then again the statement (b) must hold.
Therefore the negation of (a) implies (b), which completes the
proof of the Lemma.
\qed

\begin{lemma}
\label{diagonal}
Let $K=K(M,D)$ be a self-affine Sierpi\'nski sponge in $\R^n$ and $\mu$
the natural probability measure on it as described in 
Definition~\ref{def:K}, let 
$D_n=M^{n-1}(D)+M^{n-2}(D)+\ldots+M(D)+D$
and suppose that there exists a
$w_n\in\{-1,0,1\} \times \ldots \times \{-1,0,1\} \sm \{(0,\ldots,0)\}$
such that 
$$
D_n+w_n+M^n(\Z^n)=D_n+M^n(\Z^n).
$$

Then K is of the form $K=L\times K_0$, where $L$ is a diagonal
of a cube $[0,1]^l$, where $l\in\{1,2,\ldots,n\}$ 
and $K_0$ is a smaller dimensional self-affine Sierpi\'nski sponge.
\end{lemma}

\bp
Since every condition is invariant under any autoisometry
of the cube $[0,1]^n$ and by such a transformation we can map
$w_n$ to a vector of the form $(1,\ldots,1,0,\ldots,0)$ we
can suppose that 
$$
w_n=(\underbrace{ 1,\ldots,1 }_{l_n},0,\ldots,0),
\qquad \textrm{ where } l_n\in\{1,2,\ldots,n\}.
$$

Now we can apply Lemma~\ref{discrete} for $i=n$, $l=l_n$.
If statement (b) of Lemma~\ref{discrete} holds then 
let $l_{n-1}=l'$ and apply the lemma again for $i=n-1$, $l=l_{n-1}$.
If (b) holds again then we continue. 
Since $n\ge l_n> l_{n-1} > l_{n-2} >\ldots \ge 1$ we cannot
repeat this for more than $n-1$ times, hence for some
$1\le i \le n$ (a) of Lemma~\ref{discrete} must hold when we
apply the lemma for $i, l=l_i$.
This way we get $i,l\in\{1,\ldots,n\}$ such that
(\ref{eq:112}) and (a) of Lemma~\ref{discrete} hold.

It is easy to see that (\ref{eq:112}) implies that
$$
D+(\underbrace{ 1,\ldots,1 }_{l},0,\ldots,0)+M(Z^n)=
D+M(Z^n)
$$
and also that this and (a) of Lemma~\ref{discrete}
implies that $D$ must be of the form
$$
D=\{(\underbrace{ a,\ldots,a }_l)\ :\ a\in\{0,1,\ldots,m_1-1\}\}
\times D',
$$
where 
$
D'\su\{0,1,\ldots,m_{l+1}-1\}\times\ldots\times\{0,1,\ldots,m_{n}-1\}
$
and $m_1=\ldots=m_l$.
Then $K=K(M,D)$ must be exactly of the claimed form,
which completes the proof.
\qed

Now we are ready to characterize those self-affine sponges for which
$\mu(K\cap(K+t))$ can be positive for ``irregular'' translations.

\begin{theorem}
\label{structure}
Let $K=K(M,D)$ be a self-affine Sierpi\'nski sponge in $\R^n$ and $\mu$
the natural probability measure on it as described in 
Definition~\ref{def:K} and let $t\in\R^n$.

Then $\m{K\cap(K+t)}=0$ holds except in the following two 
trivial exceptional cases:


\noindent
(i) There exists two elementary parts $S_1$ and $S_2$ of $K$
such that $S_2=S_1+t$.

\noindent
(ii) $K$ is of the form $K=L\times K_0$, where $L$ is a diagonal
of a cube $[0,1]^l$, where $l\in\{1,2,\ldots,n\}$ 
and $K_0$ is a smaller dimensional self-affine Sierpi\'nski sponge.
\end{theorem}

\bp
If $\norma{M^k t}=0$ for some $k\in\{0,1,2,\ldots\}$
then for any two $k$-th generation elementary parts $S_1$ and $S_2$ 
of $K$,  
$S_2$ and $S_1+t$ are either identical or 
$\m{(S_1+t)\cap S_2}=0$.
Therefore in this case either (i) or $\m{K\cap (K+t)}=0$ holds,
thus we can suppose that 
$\norma{M^k t}>0$ for every $k=0,1,2,\ldots$
and $\m{K\cap (K+t)}>0$.

Let $D_i=M^{i-1}(D)+M^{i-2}(D)+\ldots+M(D)+D$. 
Notice that, by definition, $K(M,D)=K(M^i,D_i)$ for any $i\in\N$.
Therefore we can apply Proposition~\ref{prop} to $(M^n,D_n)$
to obtain $w\in\{-1,0,1\}^n\sm \{(0,\ldots,0)\}$
such that 
$$
D_n+w_n+M^n(\Z^n)=D_n+M^n(\Z^n).
$$
Then we can apply Lemma~\ref{diagonal} to get that 
$K=K(M,D)$ must be exactly of the form as in (ii) of
Theorem~\ref{structure}, which completes the proof.
\qed

\br
 Clearly, case (i) holds if and only if $t$ is
of the form
$\sum_{j=1}^k M^{-j}(\al_j-\be_j)$, 
where $k\in\{0,1,2,\ldots\}$ and $\al_1,\be_1,\dots,\al_k,\be_k\in D$.
\er

\br It follows from the proof that in the coordinates of
$L$ every $m_i$ must be the same hence in case (ii) we must
have $l=1$ if $m_1,\ldots,m_n$ are all distinct.

In particular, if $n=1$ then (ii) means $K=[0,1]$.
\er

The following statement is the analogue of Theorem~\ref{pos}.
 
\bcoro
\label{interior}
Let $K\su\R^n$ $(n\in\N)$ be a self-affine Sierpi\'nski sponge 
and $\mu$ the natural 
probability measure on it (as described in Definition~\ref{def:K}) and
let $t\in\R^n$.

The set $K\cap (K+t)$ has positive $\mu$-measure if and only if
it has non-empty interior (relative) in $K$.
\ecoro

\bp
If $K\cap (K+t)$ has non-empty interior in $K$ then clearly 
$\m{K\cap (K+t)}>0$.

We shall prove the converse by induction. Assume that the
converse is true for any smaller dimensional self-affine Sierpi\'nski sponge.
Suppose that $\m{K\cap (K+t)}>0$ and apply Theorem~\ref{structure}.
If (i) of Theorem~\ref{structure} holds then clearly 
$K\cap (K+t)$ has non-empty interior in $K$, 
so we can suppose that (ii) holds:
$K=L\times K_0$, $L$ is a diagonal of $[0,1]^l$ and $K_0$ is a
smaller dimensional self-affine Sierpi\'nski sponge. 
Then $\mu=c\la \times \mu_0$, where $1/c$ is the length
of $L$ (that is, $c=1/\sqrt l$), $\la$ is the (one-dimensional)
Lebesgue measure on $L$ and $\mu_0$ is the natural probability
measure on $K_0$.

Let $t_\al=(t_1,\ldots,t_l)$ and $t_\be=(t_{l+1},\ldots,t_n)$ and
we suppose that the coordinates of $L$ are the first $l$ coordinates.
Then
$$
K\cap(K+t)= (L\times K_0)\cap ((L+t_\al)\times (K_0+t_\be)) =
(L\cap (L+t_\al))\times (K_0\cap (K_0+t_\be)).
$$
Therefore we have 
$$
0<\m{K\cap(K+T)}=c\la\big(L\cap(L+t_\al)\big)\cdot \mu_0\big(K_0\cap(K_0+t_{\be})\big)
$$
and so $\la\big(L\cap(L+t_\al)\big)>0$ and $\mu_0\big(K_0\cap(K_0+t_{\be})\big)>0$.
This implies that $L\cap(L+t_\al)$ has non-empty interior in $L$
and, by our assumption, $K_0\cap(K_0+t_{\be})$ has non-empty interior
in $K_0$. 
Thus $K\cap(K+t)=(L\cap (L+t_\al))\times (K_0\cap (K_0+t_\be))$
has non-empty interior in $K=L\times K_0$.
\qed

For getting the analogue of Theorem~\ref{andras}
we need one more lemma.

\begin{propo}
\label{propA}
Let $K=K(M,D)$ and $\mu$ be like in Definition~\ref{def:K},
and let $0\neq t\in\R^n$ be such that $\m{K\cap (K+t)}>1-\frac1{r^2}$.

Then there exists a 
$$
w\in\{-1,0,1\} \times \ldots \times \{-1,0,1\} \sm \{(0,\ldots,0)\}
$$
such that $D+w=D$ modulo $(m_1,\ldots,m_n)$; that is,
$$
D+w+M(\Z^n)=D+M(\Z^n)
$$
\end{propo}

\bp
By Proposition~\ref{prop} we are done if 
$\norm M^k t \norm >0$ for every $k=1,2,\ldots$.
Thus we can suppose that this is not the case and choose a
minimal $k\in\{1,2,\ldots\}$ such that $\norm M^k t\norm =0$.
Then, letting $u=M^k t$, we have $u\in\Z^n\sm M(\Z^n)$.

Let 
$$
D_k=M^{k-1}(D)+M^{k-2}(D)+\ldots+M(D)+D,
$$
and define the measure $\mu_k$ so that $\mu_k(M^k A)=r^k \mu(A)$
for any Borel set $A\su K$. 
Then by definition we have $M^k K=K+D_k$, and for each $d\in D_k$
we have $\mu_k(K+d)=1$. 
Using the above facts and definitions and the condition 
$\m{K\cap (K+t)}>1-\frac1{r^2}$,
we get
\begin{multline*}
r^{k-2}(r^2-1)=r^k\Big(1-\frac1{r^2}\Big)<r^k\m{K\cap(K+t)}=
\mu_k\big(M^k K\cap (M^k K+ M^k t)\big) \\
=\mu_k\big((K+D_k)\cap (K+D_k+u)\big)=
\#(D_k\cap(D_k+u)),
\end{multline*}
where $\#(.)$ denotes the number of the elements of a set.

Then by the pigeonhole principle there exists an 
$e\in M^{k-1}(D)+M^{k-2}(D)+\dots+M^2(D)\su M^2(\Z^n)$ such that
$e+M(D)+D \su D_k+u$.
This implies that $M(D)+D+M^2(\Z^n)\su D_k+u+M^2(\Z^n)=M(D)+D+u+M^2(\Z^n)$.
Similarly, we can prove that $M(D)+D+u+M^2(\Z^n)\su M(D)+D+M^2(\Z^n)$.
Therefore we get
\begin{equation}\label{D2}
M(D)+D+u+M^2(\Z^n)=M(D)+D+M^2(\Z^n).
\end{equation}
In particular, we have
$D+u+M(\Z^n)=D+M(\Z^n).$

Then, starting from arbitrary $f_0\in D$ we can get a sequence
$(f_i)\su D$ so that 
\begin{equation}\label{sequence}
f_i+ u+M(\Z^n)= f_{i+1}+ M(\Z^n) \qquad 
(i=0,1,2,\ldots).
\end{equation}

Since $u\not\in M(\Z^n)$ we have $f_i\neq f_{i+1}$ for each $i$.
This and the fact that the sequence $(f_i)$ is contained in a finite
set imply that there must be a $j\in\N$ such that 
$f_{j+1}-f_j \neq f_j-f_{j-1}$.

Let $e\in D$ be arbitrary. Applying (\ref{D2}) and (\ref{sequence})
we get that there exist $e', e''\in D$ such that
$$
Me'+f_{j-1}+u+M^2(\Z^n)=Me+f_j+M^2(\Z^n)
$$
and
$$
Me'+f_{j}+u+M^2(\Z^n)=Me''+f_{j+1}+M^2(\Z^n),
$$
which implies
$$
(f_j-f_{j-1})-(f_{j+1}-f_j)=M(e''-e)+M^2(\Z^n).
$$

Thus there exists a $w\in\Z^n$ such that
\begin{equation}\label{Mw}
Mw=(f_j-f_{j-1})-(f_{j+1}-f_j)=M(e''-e)+M^2(\Z^n).
\end{equation}

Since $e, e'', f_{j-1},f_j, f_{j+1}\in D\su
\{0,1,\ldots,m_1-1\} \times \ldots \times \{0,1,\ldots,m_n-1\}$,
(\ref{Mw}) implies that 
$$
e+w+M(\Z^n)=e''+M(\Z^n)
$$
and
$$
w\in\{-1,0,1\} \times \ldots \times \{-1,0,1\} \sm \{(0,\ldots,0)\}.
$$
Since $e\in D$ was arbitrary, $e''\in D$ and $w$ does not depend on $e$ 
we get that
$$
D+w+M(\Z^n)=D+M(\Z^n),
$$
which completes the proof.
\qed

\bt
\label{instabilsponge}
Let $K=K(M,D)$ be a self-affine Sierpi\'nski sponge in $\R^n$ and $\mu$
the natural probability measure on it as described in 
Definition~\ref{def:K} and let $t\in\R^n$.

Then $\m{K\cap(K+t)}\le 1-\frac{1}{r^2}$ holds 
(where $r$ denotes the number of elements in the pattern $D$)
except in the following two 
trivial exceptional cases:

\noindent
(i) $t=0$.

\noindent
(ii) $K$ is of the form $K=L\times K_0$, where $L$ is a diagonal
of a cube $[0,1]^l$, where $l\in\{1,2,\ldots,n\}$ 
and $K_0$ is a smaller dimensional self-affine Sierpi\'nski sponge.
\et

\bp
Suppose that $t\neq 0$ and $\mu(K\cap(K+t))> 1-\frac{1}{r^2}$.
For $D_n=M^{n-1}(D)+M^{n-2}(D)+\ldots+M(D)+D$, 
by definition, $K(M,D)=K(M^n,D_n)$.
Therefore we can apply Proposition~\ref{propA} to $(M^n,D_n)$
to obtain $w_n\in\{-1,0,1\}^n\sm \{(0,\ldots,0)\}$
such that 
$$
D_n+w_n+M^n(\Z^n)=D_n+M^n(\Z^n).
$$
Then Lemma~\ref{diagonal} completes the proof.
\qed


\section{Translation invariant measures
for self-affine Sierpi\'nski sponges}
\label{spongemeasure}

As an easy application of Theorem~\ref{structure} (and Lemma~\ref{extendspec})
we get the following.

\bt
\label{main}
For any self-affine Sierpi\'nski sponge $K\su\Rn$ ($n\in\N$) 
there exists a translation invariant Borel measure $\nu$ on $\Rn$
such that $\nu(K)=1$. 
\et

\bp
Let $\mu$ be the natural probability Borel measure on $K$ 
(see Definition~\ref{def:K}). We shall prove by induction that
$\mu$ can be extended to $\Rn$ as a translation invariant Borel measure.
Assume that this is true for any smaller dimensional 
self-affine Sierpi\'nski sponge.

First suppose that $K$ is of the form $K=K_1\times K_2$, where 
$K_1$ and $K_2$ are smaller dimensional 
self-affine Sierpi\'nski sponges. 
Then $\mu=\mu_1\times\mu_2$, where $\mu_1$ and $\mu_2$ are the
natural probability Borel measures on $K_1$ and $K_2$, respectively.
Then, by our assumption, $\mu_1$ and $\mu_2$ has translation invariant
extensions $\mut_1$ and $\mut_2$ and then one can easily check
that $\mut=\mut_1\times\mut_2$ is a translation invariant Borel measure
on $\R^n$ and an extension of $\mu$.

If $K$ is not of the form $K=K_1\times K_2$ then we shall check that
condition (\ref{inv}) of Lemma~\ref{extendspec} is satisfied, so 
then Lemma~\ref{extendspec} will complete the proof. 
Fix $B\su K$ and $t\in\R^n$ such that $B+t\su K$. 
Then $B\su K\cap (K-t)$ and $B+t\su K\cap(K+t)$, so we have 
$\mu(B)=0=\mu(B+t)$ unless
\begin{equation}
  \label{e:pos}
  \m{K\cap(K+t)}>0 \quad \textrm{or} \quad \m{K\cap(K-t)}>0
\end{equation}
By Theorem~\ref{structure} and since case (ii) of Theorem~\ref{structure}
is already excluded, (\ref{e:pos}) implies (i) of Theorem~\ref{structure}.
On the other hand, if (i) of Theorem~\ref{structure} holds then the 
translation by $t$ maps elementary parts of $B$ to elementary parts
of $B+t$ and then the condition (\ref{inv}) clearly holds.

Since we checked all cases, the proof is complete.
\qed

We also show a more direct proof for the above theorem, which does not use 
Theorem~\ref{structure} and which works for a slightly larger class
of self-affine sets.

\bt\label{direct}
Let $\phi$ be a contractive affine map, $t_1,\ldots,t_r\in\R^n$ and
$K\su\R^n$ the compact self-affine set such that 
$K=\cup_{i=1}^r \phi(K)+t_i$.
Suppose that the standard natural probability measure on $K$ has the property
that 
\begin{equation}\label{sepcond}
\mu\Big(K\cap\Big(\big((\phi(K)+t_i)\cap(\phi(K)+t_j)\big)+u\Big)\Big)=0 \quad
(\forall\ 1\le i<j\le r,\ u\in\R^n).
\end{equation}
\begin{itemize}
\item[(a)] Then for any $t\in\R^n$ and elementary part $S$ of $K$ we have
$$
\m{K\cap(S+t)}\le \mu(S).
$$
\item[(b)] There exists a translation invariant Borel measure $\nu$ on $\R^n$
such that $\nu(K)=1$. In fact, $\nu$ is an extension of $\mu$.
\end{itemize}
\et

\bp
First we prove (a). 
Suppose that $S$ is a $k$-th generation elementary part of $K$.
Then $K$ can be written as
$$
K=\cup_{j=1}^{r^k} S+h_j
$$
for some $h_1,\ldots,h_{r^k}\in\R^n$ and by (\ref{sepcond}) the
sets $S+h_j$ are pairwise almost disjoint.

Using this and that $\mu(A)=\mu(A+h_j)$ for any Borel set $A\su S$ we
get that
\begin{align}\label{KSt}
\m{K\cap(S+t)}&=
\mu\bigg(\bigcup_{j=1}^{r^k} (S+h_j)\cap(S+t)\bigg)\nonumber\\
&= \sum_{j=1}^{r^k} \M{(S+h_j)\cap(S+t)} \nonumber\\
&= \sum_{j=1}^{r^k} \M{\big(S\cap(S+t-h_j)\big)+h_j} \nonumber\\
&= \sum_{j=1}^{r^k} \M{S\cap(S+t-h_j)}.
\end{align}

Using (\ref{sepcond}) we get that for any $i\neq j$ we have
\begin{multline}\label{commutative}
\mu\Big(\big(S\cap(S+t-h_i)\big) \cap \big(S\cap(S+t-h_j)\big)\Big)\\
=\mu\Big(S\cap\Big(\big((S+h_j)\cap(S+h_i)\big)+t-h_i-h_j\Big)\Big)=0.
\end{multline}

Thus we can continue (\ref{KSt}) as
\begin{multline*}
\m{K\cap(S+t)}=\sum_{j=1}^{r^k} \m{S\cap(S+t-h_j)}
=\mu\bigg(S\cap\bigcup_{j=1}^{r^k} (S+t-h_j)\bigg)\le \mu(S),
\end{multline*}
which completes the proof of (a).

For proving (b) define 
$$
\nu(H)=\inf\Bigg\{\sum_{j=1}^\infty \mu(S_j) : 
H\su \cup_{j=1}^\infty S_j+u_j, S_j \textrm{ is an elem. part of } K,
u_j\in\R^n\Bigg\}
$$
for any $H\su\R^n$. Then $\nu$ is clearly a translation invariant outer
measure on $\R^n$.

We claim that $\nu$ is a metric outer measure; that is, 
$\nu(A\cup B)=\nu(A)+\nu(B)$ if $A,B\su\R^n$ have positive distance.
Indeed, in this case in the cover 
$A\cup B\su \cup_{j=1}^\infty S_j+u_j$ in the definition of $\nu(A\cup B)$
we can replace replace each $S_j$ by its small elementary parts
such that each small elementary part covers only at most one of $A$ and $B$.
Since this transformation does not change $\sum_{j=1}^\infty \mu(S_j)$
this implies that $\nu(A\cup B)\ge\nu(A)+\nu(B)$. 
Since $\nu$ is an outer measure we get that 
$\nu(A\cup B)=\nu(A)+\nu(B)$.

It is well known (see e. g. in \cite{Fa85}) that 
restricting a metric outer measure to the Borel sets we get a Borel
measure.

So it is enough to prove that $\nu(K)=1$. The definition of $\nu(K)$
implies that $\nu(K)\le\mu(K)=1$.

For proving $\nu(K)\ge 1$ let $K\su\cup_{j=1}^\infty S_j+u_j$
be an arbitrary cover such that each $S_j$ is an elementary part of $K$
and $u_j\in\R^n$. Then, using the already proved (a) part we get that
$$
\sum_{j=1}^\infty \mu(S_j) \ge 
\sum_{j=1}^\infty \m{K\cap(S_j+u_j)} \ge
\mu\Big(\bigcup_{j=1}^\infty (K\cap(S_j+u_j))\Big)=\mu(K),
$$
which completes the proof of (b).
\qed

Using Lemma~\ref{cosc-msc}, the above theorem has the following consequence.

\begin{cccc}\label{cosc}
Let $\egyenlet$ be a self-affine set with the convex open set condition
and suppose that $\phi_1(K),\ldots,\phi_r(K)$ are translates of each other.

Then the natural probability measure on $K$ can be extended as a 
translation invariant measure on $\Rn$.
\ep
\end{cccc}



\section{Concluding remarks}
\label{conc}

Our results might be true for much larger classes of self-similar or
self-affine sets. We have no counter-example even for the strongest 
very naive conjecture that the intersection of any
two affine copies of \emph{any} self-affine set is of positive measure
(according to any self-affine measure on one of the copies) if and
only if it contains a set which is open in both copies.

We do not even know whether this very naive conjecture holds at least
for two isometric copies of a self-affine Sierpi\'nski sponge.
(Note that if we allow only translated copies then 
Corollary~\ref{interior} provides an affirmative answer.) 
For generalizing our results about 
Sierpi\'nski sponges from translates to isometries the following 
statement could help.

\begin{conjecture}\label{spongeconj}
If $K$ is a self-affine sponge, $\mu$ is the natural probability measure
on it, $\phi$ is an isometry and $\m{K\cap \phi(K)}>0$ then 
there exists a translation $t$ such that $K\cap\phi(K)=K\cap(K+t)$.
\end{conjecture}

This conjecture and the above mentioned Corollary~\ref{interior}
would clearly imply that Corollary~\ref{interior} holds for
isometric copies of self-affine Sierpi\'nski sponges as well.
Then, in the same way as Theorem~\ref{main} is proved, we could 
get an \emph{isometry}-invariant Borel measure $\nu$ for an arbitrary 
Sierpi\'nski sponge $K$ such that $\nu(K)=1$.

For getting this stronger version of Theorem~\ref{main} the other
natural way could be a generalization of Theorem~\ref{direct} for
isometries at least for self-affine Sierpin\'ski sponges. Since 
part (b) of Theorem~\ref{direct} follows from (a) for isometries as well
it would be enough to show (a), that is, it would be enough to
show that $\m{K\cap \phi(S)}\le\mu(S)$, for any elementary piece $S$ of
any self-affine Sierpi\'nski sponge $K$ with natural measure $\mu$.
We do not know whether this last mentioned statement holds or not.

As we saw in Theorem~\ref{instabilsponge}, the instability results
are not true for arbitrary self-affine sets, not even for
self-similar sets with the open set condition: the simplest
counter-example is 
$K=C\times [0,1]$, where $C$ denotes
the classical triadic Cantor set. Then $K$ is self-similar (with six
similitudes of ratio $1/3$), the open set condition clearly holds and if  
$\mu$ is the evenly distributed self-similar measure on $K$ (that is,
$p_1=\ldots=p_6$) then $\m{K\cap(K+(0,\eps)}=1-\eps$.
The instability results might be true for totally disconnected 
(which means that each connected component is a singleton) self-affine sets.

In the definition of self-affine sets we allowed only contractive
affine maps. If we allowed non-contractive affine maps as well 
then the above $K=C\times [0,1]$ set would be a self-affine
set (with two affine maps) with the strong separation condition,
so it would be a counter-example for both theorems 
(Theorem~\ref{kovesarki} and Theorem~\ref{selfaffine})
about self-affine sets.

We do not know whether the analogues of Theorem~\ref{andras},
Theorem~\ref{pos} and Corollary~\ref{countable} hold for self-affine
sets with the strong separation condition. Although Theorem~\ref{selfaffine}
says that for self-affine sets and isometries the analogue of 
Theorem~\ref{andras} holds, and Theorem~\ref{pos} was proved from
Theorem~\ref{andras}, we cannot get the same way that
for  self-affine sets and at least for isometries the analogue 
of Theorem~\ref{pos} holds. This is because in the proof of Theorem~\ref{pos}
it was important
that the maps $\phi_1,\ldots,\phi_r$ that generated the self-similar sets
were also in the group (in this case the group of similitudes)
for which we had Theorem~\ref{andras}. In order to get any analogue
of Theorem~\ref{pos} for self-affine sets in the same way we need
to prove a self-affine analogue of Theorem~\ref{andras} for a group
of transformation containing the affine maps $\phi_1,\ldots,\phi_r$ that 
generates the self-affine set.

From a positive answer for the following question we could get fairly easily
that the self-affine analogue of Theorem~\ref{andras} holds at least
for affine maps from any compact subset of the space of affine maps.
Then, if we could also show that we can assume that the affine maps
are from a compact set (as in Proposition~\ref{korlatos} for similitudes) 
then we would get that all the main results of Section~\ref{similar}
 also hold for self-affine
sets and affine maps as well.

\begin{question}\label{relnyilt}
Let $K\su\Rd$ be a self-affine set satisfying the strong separation condition
and let $f$ be an affine map such that $f(K)\su K$.
Does this imply that $f(K)$ is a relative open set in $K$?
\end{question}

Note that for $f(K)$ being a relative open set in $K$ means that it
is the union of countably many pairwise disjoint elementary pieces of $K$,
and since $f(K)$ is compact this means that
$f(K)$ is a finite union of elementary pieces of $K$.

A positive answer at least for the following self-similar
special case of the above question could make the
proof of Theorem~\ref{andras} simpler. However, we cannot answer this
question even for $d=1$.

\begin{question}\label{relnyilt2}
Let $K\su\Rd$ be a self-similar set satisfying the strong separation condition
and let $f$ be a similitude such that $f(K)\su K$.
Does this imply that $f(K)$ is a relative open set in $K$
(or in other words $f(K)$ is a finite union of elementary pieces of $K$)? 
\end{question}

Note that in Section~\ref{pelda} we saw that self-similar set
(even in $\R$) may contain similar copies of itself in non-trivial ways.

\end{document}

%% file: rajz.tex
\put(-.500000,0.090000){\line(1,0){1.000000}} \put(-.500000,0.085000){\line(0,1){.01}} 
\put(0.500000,0.085000){\line(0,1){.01}} 
\put(-.500000,0.000000){\line(1,0){0.150000}} \put(-.500000,-.005000){\line(0,1){.01}} 
\put(-.350000,-.005000){\line(0,1){.01}} 
\put(-.303846,0.000000){\line(1,0){0.150000}} \put(-.303846,-.005000){\line(0,1){.01}} 
\put(-.153846,-.005000){\line(0,1){.01}} 
\put(0.153846,0.000000){\line(1,0){0.150000}} \put(0.153846,-.005000){\line(0,1){.01}} 
\put(0.303846,-.005000){\line(0,1){.01}} 
\put(0.350000,0.000000){\line(1,0){0.150000}} \put(0.350000,-.005000){\line(0,1){.01}} 
\put(0.500000,-.005000){\line(0,1){.01}} 
\put(-.500000,-.090000){\line(1,0){0.022500}} \put(-.500000,-.095000){\line(0,1){.01}} 
\put(-.477500,-.095000){\line(0,1){.01}} 
\put(-.470577,-.090000){\line(1,0){0.022500}} \put(-.470577,-.095000){\line(0,1){.01}} 
\put(-.448077,-.095000){\line(0,1){.01}} 
\put(-.401923,-.090000){\line(1,0){0.022500}} \put(-.401923,-.095000){\line(0,1){.01}} 
\put(-.379423,-.095000){\line(0,1){.01}} 
\put(-.372500,-.090000){\line(1,0){0.022500}} \put(-.372500,-.095000){\line(0,1){.01}} 
\put(-.350000,-.095000){\line(0,1){.01}} 
\put(-.303846,-.090000){\line(1,0){0.022500}} \put(-.303846,-.095000){\line(0,1){.01}} 
\put(-.281346,-.095000){\line(0,1){.01}} 
\put(-.274423,-.090000){\line(1,0){0.022500}} \put(-.274423,-.095000){\line(0,1){.01}} 
\put(-.251923,-.095000){\line(0,1){.01}} 
\put(-.205769,-.090000){\line(1,0){0.022500}} \put(-.205769,-.095000){\line(0,1){.01}} 
\put(-.183269,-.095000){\line(0,1){.01}} 
\put(-.176346,-.090000){\line(1,0){0.022500}} \put(-.176346,-.095000){\line(0,1){.01}} 
\put(-.153846,-.095000){\line(0,1){.01}} 
\put(0.153846,-.090000){\line(1,0){0.022500}} \put(0.153846,-.095000){\line(0,1){.01}} 
\put(0.176346,-.095000){\line(0,1){.01}} 
\put(0.183269,-.090000){\line(1,0){0.022500}} \put(0.183269,-.095000){\line(0,1){.01}} 
\put(0.205769,-.095000){\line(0,1){.01}} 
\put(0.251923,-.090000){\line(1,0){0.022500}} \put(0.251923,-.095000){\line(0,1){.01}} 
\put(0.274423,-.095000){\line(0,1){.01}} 
\put(0.281346,-.090000){\line(1,0){0.022500}} \put(0.281346,-.095000){\line(0,1){.01}} 
\put(0.303846,-.095000){\line(0,1){.01}} 
\put(0.350000,-.090000){\line(1,0){0.022500}} \put(0.350000,-.095000){\line(0,1){.01}} 
\put(0.372500,-.095000){\line(0,1){.01}} 
\put(0.379423,-.090000){\line(1,0){0.022500}} \put(0.379423,-.095000){\line(0,1){.01}} 
\put(0.401923,-.095000){\line(0,1){.01}} 
\put(0.448077,-.090000){\line(1,0){0.022500}} \put(0.448077,-.095000){\line(0,1){.01}} 
\put(0.470577,-.095000){\line(0,1){.01}} 
\put(0.477500,-.090000){\line(1,0){0.022500}} \put(0.477500,-.095000){\line(0,1){.01}} 
\put(0.500000,-.095000){\line(0,1){.01}} 
\put(-.500000,-.180000){\line(1,0){0.022500}} \put(-.500000,-.185000){\line(0,1){.01}} 
\put(-.477500,-.185000){\line(0,1){.01}} 
\put(-.470577,-.180000){\line(1,0){0.022500}} \put(-.470577,-.185000){\line(0,1){.01}} 
\put(-.448077,-.185000){\line(0,1){.01}} 
\put(-.401923,-.180000){\line(1,0){0.150000}} \put(-.401923,-.185000){\line(0,1){.01}} 
\put(-.251923,-.185000){\line(0,1){.01}} 
\put(-.205769,-.180000){\line(1,0){0.022500}} \put(-.205769,-.185000){\line(0,1){.01}} 
\put(-.183269,-.185000){\line(0,1){.01}} 
\put(-.176346,-.180000){\line(1,0){0.022500}} \put(-.176346,-.185000){\line(0,1){.01}} 
\put(-.153846,-.185000){\line(0,1){.01}} 
\put(0.153846,-.180000){\line(1,0){0.022500}} \put(0.153846,-.185000){\line(0,1){.01}} 
\put(0.176346,-.185000){\line(0,1){.01}} 
\put(0.183269,-.180000){\line(1,0){0.022500}} \put(0.183269,-.185000){\line(0,1){.01}} 
\put(0.205769,-.185000){\line(0,1){.01}} 
\put(0.251923,-.180000){\line(1,0){0.150000}} \put(0.251923,-.185000){\line(0,1){.01}} 
\put(0.401923,-.185000){\line(0,1){.01}} 
\put(0.448077,-.180000){\line(1,0){0.022500}} \put(0.448077,-.185000){\line(0,1){.01}} 
\put(0.470577,-.185000){\line(0,1){.01}} 
\put(0.477500,-.180000){\line(1,0){0.022500}} \put(0.477500,-.185000){\line(0,1){.01}} 
\put(0.500000,-.185000){\line(0,1){.01}} 
\put(-.500000,-.270000){\line(1,0){0.000506}} 
\put(-.499338,-.270000){\line(1,0){0.000506}} 
\put(-.497793,-.270000){\line(1,0){0.000506}} 
\put(-.497131,-.270000){\line(1,0){0.000506}} 
\put(-.495587,-.270000){\line(1,0){0.000506}} 
\put(-.494925,-.270000){\line(1,0){0.000506}} 
\put(-.493380,-.270000){\line(1,0){0.000506}} 
\put(-.492718,-.270000){\line(1,0){0.000506}} 
\put(-.485288,-.270000){\line(1,0){0.000506}} 
\put(-.484626,-.270000){\line(1,0){0.000506}} 
\put(-.483082,-.270000){\line(1,0){0.000506}} 
\put(-.482420,-.270000){\line(1,0){0.000506}} 
\put(-.480875,-.270000){\line(1,0){0.000506}} 
\put(-.480213,-.270000){\line(1,0){0.000506}} 
\put(-.478668,-.270000){\line(1,0){0.000506}} 
\put(-.478006,-.270000){\line(1,0){0.000506}} 
\put(-.470577,-.270000){\line(1,0){0.000506}} 
\put(-.469915,-.270000){\line(1,0){0.000506}} 
\put(-.468370,-.270000){\line(1,0){0.000506}} 
\put(-.467708,-.270000){\line(1,0){0.000506}} 
\put(-.466163,-.270000){\line(1,0){0.000506}} 
\put(-.465501,-.270000){\line(1,0){0.000506}} 
\put(-.463957,-.270000){\line(1,0){0.000506}} 
\put(-.463295,-.270000){\line(1,0){0.000506}} 
\put(-.455865,-.270000){\line(1,0){0.000506}} 
\put(-.455203,-.270000){\line(1,0){0.000506}} 
\put(-.453659,-.270000){\line(1,0){0.000506}} 
\put(-.452997,-.270000){\line(1,0){0.000506}} 
\put(-.451452,-.270000){\line(1,0){0.000506}} 
\put(-.450790,-.270000){\line(1,0){0.000506}} 
\put(-.449245,-.270000){\line(1,0){0.000506}} 
\put(-.448583,-.270000){\line(1,0){0.000506}} 
\put(-.401923,-.270000){\line(1,0){0.000506}} 
\put(-.401261,-.270000){\line(1,0){0.000506}} 
\put(-.399716,-.270000){\line(1,0){0.000506}} 
\put(-.399054,-.270000){\line(1,0){0.000506}} 
\put(-.397510,-.270000){\line(1,0){0.000506}} 
\put(-.396848,-.270000){\line(1,0){0.000506}} 
\put(-.395303,-.270000){\line(1,0){0.000506}} 
\put(-.394641,-.270000){\line(1,0){0.000506}} 
\put(-.387212,-.270000){\line(1,0){0.000506}} 
\put(-.386550,-.270000){\line(1,0){0.000506}} 
\put(-.385005,-.270000){\line(1,0){0.000506}} 
\put(-.384343,-.270000){\line(1,0){0.000506}} 
\put(-.382798,-.270000){\line(1,0){0.000506}} 
\put(-.382136,-.270000){\line(1,0){0.000506}} 
\put(-.380591,-.270000){\line(1,0){0.000506}} 
\put(-.379929,-.270000){\line(1,0){0.000506}} 
\put(-.372500,-.270000){\line(1,0){0.000506}} 
\put(-.371838,-.270000){\line(1,0){0.000506}} 
\put(-.370293,-.270000){\line(1,0){0.000506}} 
\put(-.369631,-.270000){\line(1,0){0.000506}} 
\put(-.368087,-.270000){\line(1,0){0.000506}} 
\put(-.367425,-.270000){\line(1,0){0.000506}} 
\put(-.365880,-.270000){\line(1,0){0.000506}} 
\put(-.365218,-.270000){\line(1,0){0.000506}} 
\put(-.357788,-.270000){\line(1,0){0.000506}} 
\put(-.357126,-.270000){\line(1,0){0.000506}} 
\put(-.355582,-.270000){\line(1,0){0.000506}} 
\put(-.354920,-.270000){\line(1,0){0.000506}} 
\put(-.353375,-.270000){\line(1,0){0.000506}} 
\put(-.352713,-.270000){\line(1,0){0.000506}} 
\put(-.351168,-.270000){\line(1,0){0.000506}} 
\put(-.350506,-.270000){\line(1,0){0.000506}} 
\put(-.303846,-.270000){\line(1,0){0.000506}} 
\put(-.303184,-.270000){\line(1,0){0.000506}} 
\put(-.301639,-.270000){\line(1,0){0.000506}} 
\put(-.300977,-.270000){\line(1,0){0.000506}} 
\put(-.299433,-.270000){\line(1,0){0.000506}} 
\put(-.298771,-.270000){\line(1,0){0.000506}} 
\put(-.297226,-.270000){\line(1,0){0.000506}} 
\put(-.296564,-.270000){\line(1,0){0.000506}} 
\put(-.289135,-.270000){\line(1,0){0.000506}} 
\put(-.288473,-.270000){\line(1,0){0.000506}} 
\put(-.286928,-.270000){\line(1,0){0.000506}} 
\put(-.286266,-.270000){\line(1,0){0.000506}} 
\put(-.284721,-.270000){\line(1,0){0.000506}} 
\put(-.284059,-.270000){\line(1,0){0.000506}} 
\put(-.282514,-.270000){\line(1,0){0.000506}} 
\put(-.281852,-.270000){\line(1,0){0.000506}} 
\put(-.274423,-.270000){\line(1,0){0.000506}} 
\put(-.273761,-.270000){\line(1,0){0.000506}} 
\put(-.272216,-.270000){\line(1,0){0.000506}} 
\put(-.271554,-.270000){\line(1,0){0.000506}} 
\put(-.270010,-.270000){\line(1,0){0.000506}} 
\put(-.269348,-.270000){\line(1,0){0.000506}} 
\put(-.267803,-.270000){\line(1,0){0.000506}} 
\put(-.267141,-.270000){\line(1,0){0.000506}} 
\put(-.259712,-.270000){\line(1,0){0.000506}} 
\put(-.259050,-.270000){\line(1,0){0.000506}} 
\put(-.257505,-.270000){\line(1,0){0.000506}} 
\put(-.256843,-.270000){\line(1,0){0.000506}} 
\put(-.255298,-.270000){\line(1,0){0.000506}} 
\put(-.254636,-.270000){\line(1,0){0.000506}} 
\put(-.253091,-.270000){\line(1,0){0.000506}} 
\put(-.252429,-.270000){\line(1,0){0.000506}} 
\put(-.205769,-.270000){\line(1,0){0.000506}} 
\put(-.205107,-.270000){\line(1,0){0.000506}} 
\put(-.203562,-.270000){\line(1,0){0.000506}} 
\put(-.202900,-.270000){\line(1,0){0.000506}} 
\put(-.201356,-.270000){\line(1,0){0.000506}} 
\put(-.200694,-.270000){\line(1,0){0.000506}} 
\put(-.199149,-.270000){\line(1,0){0.000506}} 
\put(-.198487,-.270000){\line(1,0){0.000506}} 
\put(-.191058,-.270000){\line(1,0){0.000506}} 
\put(-.190396,-.270000){\line(1,0){0.000506}} 
\put(-.188851,-.270000){\line(1,0){0.000506}} 
\put(-.188189,-.270000){\line(1,0){0.000506}} 
\put(-.186644,-.270000){\line(1,0){0.000506}} 
\put(-.185982,-.270000){\line(1,0){0.000506}} 
\put(-.184437,-.270000){\line(1,0){0.000506}} 
\put(-.183775,-.270000){\line(1,0){0.000506}} 
\put(-.176346,-.270000){\line(1,0){0.000506}} 
\put(-.175684,-.270000){\line(1,0){0.000506}} 
\put(-.174139,-.270000){\line(1,0){0.000506}} 
\put(-.173477,-.270000){\line(1,0){0.000506}} 
\put(-.171933,-.270000){\line(1,0){0.000506}} 
\put(-.171271,-.270000){\line(1,0){0.000506}} 
\put(-.169726,-.270000){\line(1,0){0.000506}} 
\put(-.169064,-.270000){\line(1,0){0.000506}} 
\put(-.161635,-.270000){\line(1,0){0.000506}} 
\put(-.160973,-.270000){\line(1,0){0.000506}} 
\put(-.159428,-.270000){\line(1,0){0.000506}} 
\put(-.158766,-.270000){\line(1,0){0.000506}} 
\put(-.157221,-.270000){\line(1,0){0.000506}} 
\put(-.156559,-.270000){\line(1,0){0.000506}} 
\put(-.155014,-.270000){\line(1,0){0.000506}} 
\put(-.154352,-.270000){\line(1,0){0.000506}} 
\put(0.153846,-.270000){\line(1,0){0.000506}} 
\put(0.154508,-.270000){\line(1,0){0.000506}} 
\put(0.156053,-.270000){\line(1,0){0.000506}} 
\put(0.156715,-.270000){\line(1,0){0.000506}} 
\put(0.158260,-.270000){\line(1,0){0.000506}} 
\put(0.158922,-.270000){\line(1,0){0.000506}} 
\put(0.160466,-.270000){\line(1,0){0.000506}} 
\put(0.161128,-.270000){\line(1,0){0.000506}} 
\put(0.168558,-.270000){\line(1,0){0.000506}} 
\put(0.169220,-.270000){\line(1,0){0.000506}} 
\put(0.170764,-.270000){\line(1,0){0.000506}} 
\put(0.171426,-.270000){\line(1,0){0.000506}} 
\put(0.172971,-.270000){\line(1,0){0.000506}} 
\put(0.173633,-.270000){\line(1,0){0.000506}} 
\put(0.175178,-.270000){\line(1,0){0.000506}} 
\put(0.175840,-.270000){\line(1,0){0.000506}} 
\put(0.183269,-.270000){\line(1,0){0.000506}} 
\put(0.183931,-.270000){\line(1,0){0.000506}} 
\put(0.185476,-.270000){\line(1,0){0.000506}} 
\put(0.186138,-.270000){\line(1,0){0.000506}} 
\put(0.187683,-.270000){\line(1,0){0.000506}} 
\put(0.188345,-.270000){\line(1,0){0.000506}} 
\put(0.189889,-.270000){\line(1,0){0.000506}} 
\put(0.190551,-.270000){\line(1,0){0.000506}} 
\put(0.197981,-.270000){\line(1,0){0.000506}} 
\put(0.198643,-.270000){\line(1,0){0.000506}} 
\put(0.200188,-.270000){\line(1,0){0.000506}} 
\put(0.200850,-.270000){\line(1,0){0.000506}} 
\put(0.202394,-.270000){\line(1,0){0.000506}} 
\put(0.203056,-.270000){\line(1,0){0.000506}} 
\put(0.204601,-.270000){\line(1,0){0.000506}} 
\put(0.205263,-.270000){\line(1,0){0.000506}} 
\put(0.251923,-.270000){\line(1,0){0.000506}} 
\put(0.252585,-.270000){\line(1,0){0.000506}} 
\put(0.254130,-.270000){\line(1,0){0.000506}} 
\put(0.254792,-.270000){\line(1,0){0.000506}} 
\put(0.256337,-.270000){\line(1,0){0.000506}} 
\put(0.256999,-.270000){\line(1,0){0.000506}} 
\put(0.258543,-.270000){\line(1,0){0.000506}} 
\put(0.259205,-.270000){\line(1,0){0.000506}} 
\put(0.266635,-.270000){\line(1,0){0.000506}} 
\put(0.267297,-.270000){\line(1,0){0.000506}} 
\put(0.268841,-.270000){\line(1,0){0.000506}} 
\put(0.269503,-.270000){\line(1,0){0.000506}} 
\put(0.271048,-.270000){\line(1,0){0.000506}} 
\put(0.271710,-.270000){\line(1,0){0.000506}} 
\put(0.273255,-.270000){\line(1,0){0.000506}} 
\put(0.273917,-.270000){\line(1,0){0.000506}} 
\put(0.281346,-.270000){\line(1,0){0.000506}} 
\put(0.282008,-.270000){\line(1,0){0.000506}} 
\put(0.283553,-.270000){\line(1,0){0.000506}} 
\put(0.284215,-.270000){\line(1,0){0.000506}} 
\put(0.285760,-.270000){\line(1,0){0.000506}} 
\put(0.286422,-.270000){\line(1,0){0.000506}} 
\put(0.287966,-.270000){\line(1,0){0.000506}} 
\put(0.288628,-.270000){\line(1,0){0.000506}} 
\put(0.296058,-.270000){\line(1,0){0.000506}} 
\put(0.296720,-.270000){\line(1,0){0.000506}} 
\put(0.298264,-.270000){\line(1,0){0.000506}} 
\put(0.298926,-.270000){\line(1,0){0.000506}} 
\put(0.300471,-.270000){\line(1,0){0.000506}} 
\put(0.301133,-.270000){\line(1,0){0.000506}} 
\put(0.302678,-.270000){\line(1,0){0.000506}} 
\put(0.303340,-.270000){\line(1,0){0.000506}} 
\put(0.350000,-.270000){\line(1,0){0.000506}} 
\put(0.350662,-.270000){\line(1,0){0.000506}} 
\put(0.352207,-.270000){\line(1,0){0.000506}} 
\put(0.352869,-.270000){\line(1,0){0.000506}} 
\put(0.354414,-.270000){\line(1,0){0.000506}} 
\put(0.355076,-.270000){\line(1,0){0.000506}} 
\put(0.356620,-.270000){\line(1,0){0.000506}} 
\put(0.357282,-.270000){\line(1,0){0.000506}} 
\put(0.364712,-.270000){\line(1,0){0.000506}} 
\put(0.365374,-.270000){\line(1,0){0.000506}} 
\put(0.366918,-.270000){\line(1,0){0.000506}} 
\put(0.367580,-.270000){\line(1,0){0.000506}} 
\put(0.369125,-.270000){\line(1,0){0.000506}} 
\put(0.369787,-.270000){\line(1,0){0.000506}} 
\put(0.371332,-.270000){\line(1,0){0.000506}} 
\put(0.371994,-.270000){\line(1,0){0.000506}} 
\put(0.379423,-.270000){\line(1,0){0.000506}} 
\put(0.380085,-.270000){\line(1,0){0.000506}} 
\put(0.381630,-.270000){\line(1,0){0.000506}} 
\put(0.382292,-.270000){\line(1,0){0.000506}} 
\put(0.383837,-.270000){\line(1,0){0.000506}} 
\put(0.384499,-.270000){\line(1,0){0.000506}} 
\put(0.386043,-.270000){\line(1,0){0.000506}} 
\put(0.386705,-.270000){\line(1,0){0.000506}} 
\put(0.394135,-.270000){\line(1,0){0.000506}} 
\put(0.394797,-.270000){\line(1,0){0.000506}} 
\put(0.396341,-.270000){\line(1,0){0.000506}} 
\put(0.397003,-.270000){\line(1,0){0.000506}} 
\put(0.398548,-.270000){\line(1,0){0.000506}} 
\put(0.399210,-.270000){\line(1,0){0.000506}} 
\put(0.400755,-.270000){\line(1,0){0.000506}} 
\put(0.401417,-.270000){\line(1,0){0.000506}} 
\put(0.448077,-.270000){\line(1,0){0.000506}} 
\put(0.448739,-.270000){\line(1,0){0.000506}} 
\put(0.450284,-.270000){\line(1,0){0.000506}} 
\put(0.450946,-.270000){\line(1,0){0.000506}} 
\put(0.452490,-.270000){\line(1,0){0.000506}} 
\put(0.453152,-.270000){\line(1,0){0.000506}} 
\put(0.454697,-.270000){\line(1,0){0.000506}} 
\put(0.455359,-.270000){\line(1,0){0.000506}} 
\put(0.462789,-.270000){\line(1,0){0.000506}} 
\put(0.463451,-.270000){\line(1,0){0.000506}} 
\put(0.464995,-.270000){\line(1,0){0.000506}} 
\put(0.465657,-.270000){\line(1,0){0.000506}} 
\put(0.467202,-.270000){\line(1,0){0.000506}} 
\put(0.467864,-.270000){\line(1,0){0.000506}} 
\put(0.469409,-.270000){\line(1,0){0.000506}} 
\put(0.470071,-.270000){\line(1,0){0.000506}} 
\put(0.477500,-.270000){\line(1,0){0.000506}} 
\put(0.478162,-.270000){\line(1,0){0.000506}} 
\put(0.479707,-.270000){\line(1,0){0.000506}} 
\put(0.480369,-.270000){\line(1,0){0.000506}} 
\put(0.481914,-.270000){\line(1,0){0.000506}} 
\put(0.482576,-.270000){\line(1,0){0.000506}} 
\put(0.484120,-.270000){\line(1,0){0.000506}} 
\put(0.484782,-.270000){\line(1,0){0.000506}} 
\put(0.492212,-.270000){\line(1,0){0.000506}} 
\put(0.492874,-.270000){\line(1,0){0.000506}} 
\put(0.494418,-.270000){\line(1,0){0.000506}} 
\put(0.495080,-.270000){\line(1,0){0.000506}} 
\put(0.496625,-.270000){\line(1,0){0.000506}} 
\put(0.497287,-.270000){\line(1,0){0.000506}} 
\put(0.498832,-.270000){\line(1,0){0.000506}} 
\put(0.499494,-.270000){\line(1,0){0.000506}} 

%% file: submittedetds.bbl
\begin{thebibliography}{111}
\bibitem{BG} C.~Bandt and S.~Graf, 
Self-similar sets. VII. 
A characterization of self-similar fractals with positive Hausdorff measure, 
\textsl{Proc. Amer. Math. Soc.}  \textbf{114}  (1992),  no. 4, 995--1001.

\bibitem{Be} T.~Bedford, 
\textsl{Crinkly curves, Markov partitions and box dimension 
in self-similar sets}, Ph.~D.~Thesis, University of Warwick, 1984.

\bibitem{Da} R.~O.~Davies, Sets which are null or non-sigma-finite
for every trans\-la\-tion-invariant measure,
\textsl{Mathematika} \textbf{18} (1971), 161--162.

\bibitem{EK} M.~Elekes and T.~Keleti, 
Borel sets which are null or non-$\sigma$-finite
for every translation invariant measure, 
\emph{Adv. Math.} \textbf{201} (2006), 102-115.

\bibitem{Fa}
K.~J.~Falconer,
Classes of sets with large intersection,
\emph{Mathematika} \textbf{32} (1985), no. 2, 191--205.

\bibitem{Fa85} K.~J.~Falconer, \textsl{The geometry of fractal sets.} 
Cambridge Tracts in Mathematics No. 85, Cambridge University Press, 1985. 

\bibitem{Fa97}
K.~J.~Falconer,
\emph{Techniques in Fractal Geometry},
John Wiley \& Sons, 1997.

\bibitem{FW}
D.~Feng and Y.~Wang,
\emph{On the structures of generating iterated function systems of 
Cantor sets},
preprint.

\bibitem{Fu}
H.~Furstenberg,
Intersections of Cantor sets and transversality of semigroups,
\emph{Problems in analysis 
(Sympos. Salomon Bochner, Princeton Univ., Princeton, N.J., 1969)},
pp. 41-59, Princeton Univ. Press, 1970.

\bibitem{GL} D.~Gatzouras and S.~Lalley (1992),
Hausdorff and box dimensions of certain self-affine fractals,
\textsl{Indiana University Math. J.}  \textbf{41} (1992), 533--568.

\bibitem{H}
{J.~E.~Hutchinson},
Fractals and self-similarity,
\emph{Indiana Univ. Math. J.} \textbf{30} (1981), no. 5, 713--747.

\bibitem{I}
K.~Igudesman,
Lacunary self-similar fractal sets and its application to
intersection of Cantor sets,
\emph{Lobachevskii J. Math.} \textbf{12} (2003), 41--50.

\bibitem{J}
M.~J\"arvenp\"a\"a,
Hausdorff and packing dimensions, intersection measures, and similarities,
\emph{Ann. Acad. Sci. Fenn. Math.} \textbf{24} (1999), no. 1., 165--186.

\bibitem{Ke} T.~Keleti,
A 1-dimensional subset of the reals that intersects 
each of its translates in at most a single point,
\textsl{Real Analysis Exchange} \textbf {24} (1998/99), 843--844.

\bibitem{KP} R.~Kenyon and Y.~Peres, 
 Measures of full dimension on affine-invariant sets,
\textsl{Ergodic Theory Dynamical Syst. }  \textbf{16} (1996), 307--323.

\bibitem{LW} J.~C.~Lagarias and Y.~Wang,  
Self-affine tiles in $ R\sp n$,
\textsl{Adv. Math.}  \textbf{121}  (1996),  no. 1, 21--49.

\bibitem{LX}
{W.~Li and D.~Xiao},
Intersection of translations of Cantor triadic set,
\emph{Acta Math. Sci. (English Ed.)} \textbf{19} (1999), no. 2, 214--219.

\bibitem{M}
A.~M\'ath\'e,
\emph{\"Onhasonl\'o halmazok egybev\'ag\'os\'ag-invari\'ans
  m\'ert\'ekeir\H{o}l} (On isometry invariant measures of self-similar sets)
(In Hungarian), Master Thesis, E\"otv\"os Lor\'and University, 2005 
({\tt http://www.cs.elte.hu/math/diploma/math}).

\bibitem{M1}
P.~Mattila,
On the structure of self-similar fractals,
\emph{Ann. Acad. Sci. Fenn. Ser. A I Math.} \textbf{7} (1982), no. 2., 189--195.

\bibitem{M2}
P.~Mattila,
Hausdorff dimensions and capacities of intersections of sets in $n$-space,
\emph{Acta Math.} \textbf{152} (1984), no. 1-2, 77--105.

\bibitem{M3}
P.~Mattila,
On the Hausdorff dimension and capacities of intersections,
\emph{Mathematika} \textbf{32} (1985), no. 2, 213--217.

\bibitem{Mu}
C.~McMullen, 
The Hausdorff dimension of general Sierpi\'nski carpets,
\textsl{Nagoya Math. J.} \textbf{96} (1984), 1--9.

\bibitem{Mo}
C.~G.~T.~de~A.~Moreira,
Stable intersections of Cantor sets and homoclinic bifurcations,
\emph{Ann. Inst. H. Poincar\'e Anal. Non Lin\'eaire} \textbf{13}, no. 6, 741--781.

\bibitem{MY}
C.~G.~T.~de~A.~Moreira and J-C.~Yoccoz, 
Stable intersections of regular Cantor sets 
with large Hausdorff dimensions,  
\textsl{Ann. of Math. (2)  }  \textbf{154}  (2001),  no. 1, 45--96.

\bibitem{NL} F.~Nekka and J.~Li, 
Intersection of triadic Cantor sets with their translates. I. 
Fundamental properties,  
\textsl{Chaos Solitons Fractals}  \textbf{13}  (2002),  no. 9, 1807--1817.

\bibitem{PS}
{Y.~Peres and B.~Solomyak},
Self-similar measures and intersections of Cantor sets,
\emph{Trans. Amer. Math. Soc.} \textbf{350} (1998), no. 10, 4065--4087.

\bibitem{PeP} Y.~Peres, 
The packing measure of self-affine carpets, 
\textsl{Math. Proc. Cambridge Philos. Soc.} 
\textbf{115} (1994), no.~3, 437--450.

\bibitem{PeH} Y.~Peres,
The self-affine carpets of McMullen and Bedford have infinite Hausdorff
   measure,
\textsl{Math. Proc. Cambridge Philos. Soc.} 
\textbf{116} (1994), 513--526. 

\bibitem{Ru} W.~Rudin, 
\textsl{Real and complex analysis}, 
Third Edition, McGraw-Hill Book Company, 1987.

\bibitem{S} A.~Schief,
Spearation properties for self-similar sets,
\emph{Proc. Amer. Math. Soc.} \textbf{122} (1994), no. 1, 111--115.




\end{thebibliography}
